\documentclass[12pt,oneside]{amsart}

\usepackage{hyperref}
\hypersetup{
    colorlinks=true,
    linkcolor=black,
    citecolor=magenta,
    filecolor=black,      
    urlcolor=blue,
    }

\usepackage{fullpage}
\usepackage{amssymb}
\usepackage{latexsym}
\usepackage{amsmath}
\usepackage{mathrsfs}
\usepackage[all]{xy}
\usepackage{enumerate}
\usepackage{enumitem}
\usepackage{amsthm}
\usepackage{psfrag}
\usepackage{ifthen}
\usepackage{mathdots}
\usepackage{caption}
\usepackage{subfigure}
\usepackage{standalone}
\usepackage{mathtools}
\usepackage{pgfplots}
\usepackage{comment}
\usepackage{xfrac}
\usepackage{todonotes}
\usepackage[utf8]{inputenc}
\usepackage[linesnumbered,ruled,vlined]{algorithm2e}
\usepackage[outline]{contour}
\usepackage{float}
\usepackage{tikz}
\usetikzlibrary{shapes.geometric}
\usepackage{graphicx}
\newcommand{\newbold}[1]{\bgroup\contourlength{0.01em}\contour{black}{#1}\egroup}

\newtheorem{theorem}{Theorem}[section]
\newtheorem{lemma}[theorem]{Lemma}
\newtheorem{metalemma}[theorem]{Meta Lemma}
\newtheorem{proposition}[theorem]{Proposition}
\newtheorem{corollary}[theorem]{Corollary}

\newtheorem*{theorem*}{Theorem}          

\theoremstyle{definition}
\newtheorem{definition}[theorem]{Definition}

\newtheorem{example}[theorem]{Example}

\usepackage{tabularx,environ}

\makeatletter
\newcommand{\problemtitle}[1]{\gdef\@problemtitle{#1}}
\newcommand{\probleminput}[1]{\gdef\@probleminput{#1}}
\newcommand{\problemquestion}[1]{\gdef\@problemquestion{#1}}
\NewEnviron{algproblem}{
  \problemtitle{}\probleminput{}\problemquestion{}
  \BODY
  \par\addvspace{.5\baselineskip}
  \noindent
  \begin{tabularx}{\textwidth}{@{\hspace{0em}} l X c}
    \hline\hline
    \multicolumn{2}{@{\hspace{0em}}l}{\@problemtitle} \\
    \textbf{In:} & \@probleminput \\%
    \textbf{Question:} & \@problemquestion\\
    \hline\hline
  \end{tabularx}
  \par\addvspace{.5\baselineskip}
}
\makeatother

\newcommand{\bel}[1]{\begin{equation}\label{#1}}
\newcommand{\ee}{\end{equation}}

\newcommand{\LBA}{\left\{\begin{array}}
\newcommand{\EAR}{\end{array}\right.}

\def\ove{{\mathbf{e}}}

\def\ovu{{\mathbf{u}}}

\def\NP{{\mathbf{NP}}}
\def\DP{{\mathbf{DP}}}

\def\blf{{\protect\newbold{$f$}}}
\def\blg{{\protect\newbold{$g$}}}
\def\blo{{\protect\newbold{$0$}}}

\def\one{{\mathbf{1}}}
\def\SD{{\,\triangle\,}}

\newcommand{\gpr}[2]{{\left\langle #1 \mid #2 \right\rangle}}
\newcommand{\rb}[1]{{\left( #1 \right)}}

\newcommand{\Set}[2]{\left\{\, #1 \;\middle|\; #2 \,\right\}}

\def\ME{{\mathbb{E}}}
\def\MN{{\mathbb{N}}}
\def\MZ{{\mathbb{Z}}}

\DeclareMathOperator{\DIV}{{DIV}}

\DeclareMathOperator{\size}{{size}}
\DeclareMathOperator{\supp}{{supp}}

\DeclareMathOperator{\ord}{{ord}}
\DeclareMathOperator{\num}{{\protect\newbold{num}}}
\DeclareMathOperator{\den}{{\protect\newbold{den}}}
\DeclareMathOperator{\pres}{{pres}}

\DeclareMathOperator{\Path}{{path}}
\DeclareMathOperator{\rev}{{rev}}
\DeclareMathOperator{\cont}{{cont}}
\DeclareMathOperator{\dmin}{{dmin}}
\DeclareMathOperator{\dmax}{{dmax}}
\DeclareMathOperator{\wO}{{\widetilde{O}}}

\title{One-variable equations over the lamplighter group}
\author{Alexander Ushakov and Yankun Wang}
\thanks{The authors are grateful to an anonymous reviewer whose suggestions to the divisibility problem led to substantial simplifications and improved complexity estimates in this paper.}

\address{Department of Mathematical Sciences, Stevens Institute of Technology, Hoboken NJ 07030}\email{aushakov,ywang7@stevens.edu}

\date{\today}

\begin{document}

\maketitle

\begin{abstract}
We study one-variable equations over the lamplighter group
$L_2=\MZ_2 \wr \MZ$.
While the decidability of arbitrary equations over 
$L_2$ remains open,
we prove that the Diophantine problem for single equations in one variable
is decidable.
Our approach reduces the problem to a divisibility question for families of
parametric Laurent polynomials over $\MZ_2$, whose exponents depend linearly
on an integer parameter.
To analyze this divisibility problem, we introduce a symbolic division
procedure for associated bivariate polynomials and derive explicit bounds on
the parameter from the structure of the resulting quotient and remainder.
This yields an explicit decision procedure with exponential worst-case
complexity.
On the other hand, we show that for a generic class of equations, solvability
can be decided in nearly quadratic time.
These results establish a sharp contrast between worst-case and typical
computational behavior and provide new tools for the study of equations over
wreath products.

\end{abstract}

\section{Introduction}

Let $F = F(Z)$ denote the free group on countably many generators 
$Z = \{z_i\}_{i=1}^\infty$. For a group $G$, an \emph{equation over $G$ with variables in $Z$} is an equation of the form $w = 1$, where $w \in F*G$. If $w = z_{i_1}g_1\cdots z_{i_k}g_k$, with $z_{i_j}\in Z$ and $g_j\in G$, then we refer to $\{z_{i_1},\ldots,z_{i_k}\}$ as the set of \emph{variables} and to $\{g_1,\ldots,g_k\}$ as the set of \emph{constants} (or \emph{coefficients}) of $w$.

A \emph{solution} to an equation $w(z_1,\ldots,z_k)=1$ over $G$ is a homomorphism $$\varphi\colon F*G\rightarrow G $$
such that $\varphi|_G = 1_G$ and $w \in \ker \varphi$. If $\varphi$ is a solution of $w =1$ and $g_i = \varphi(z_i)$, then we often say that $g_1,\ldots,g_k$ is a solution of $w =1$.

We assume that $G$ comes equipped with a fixed generating set $X$ 
and elements of $G$ are given as products of elements of $X$ and their inverses. 
This naturally defines the length (or size) of the equation $w =1$ as the length 
of its left-hand side word $w $.

The \emph{Diophantine problem}, $\DP_C(G)$, over a group $G$ 
for a class of equations $C$ 
is an algorithmic question to decide whether a given equation $w =1$ in $C$ 
has a solution in $G$.
It can be regarded as a group-theoretic analogue of the satisfiability problem.
More generally, one can study the Diophantine problem for systems of
equations.
In this paper, we focus on the class of one-variable equations over 
the lamplighter group $\MZ_2\wr \MZ$.

We say that $w =1$ is a \emph{one-variable equation} if $w$ involves a single
variable, for convenience we assume it is $x$, one or more times
as $x$ or as $x^{-1}$.

\begin{algproblem}
\problemtitle{\textsc{Diophantine problem for one-variable equations over $G$} $(\DP_1(G))$.}
\probleminput{A group word $w =w(X,x) \in F(x)\ast G$, where $X$ is a generating set for $G$.}
\problemquestion{Is there a group word $w(X)$ such that $w(X,w)=1$ in $G$?}
\end{algproblem}

\subsection{Equations in the lamplighter and similar groups}

Computational properties of equations in the lamplighter group, 
in wreath products, and more generally 
in finitely generated metabelian groups have attracted 
considerable attention in recent years.

A.~Miasnikov and N.~Romanovsky \cite{Myasnikov-Romanovskii:2012} proved
that the Diophantine problem for systems of equations over $\MZ\wr \MZ$ 
is undecidable. 
R.~Dong \cite{Dong:2025} proved the same result using a different argument.
However, it is unknown if the problem for single equations
over $\MZ\wr \MZ$ is decidable or not.
I.~Lysenok and A.~Ushakov proved that the Diophantine problem 
over free metabelian groups is $\NP$-complete for 
spherical equations \cite{Lysenok-Ushakov:2015}, and more generally, 
for orientable equations \cite{Lysenok-Ushakov:2021}.
The Diophantine problem 
for orientable quadratic equations over wreath products 
of finitely generated abelian groups is $\NP$-complete \cite{Ushakov-Weiers:2025}.
Also, the Diophantine problem for quadratic equations is $\NP$-complete
over metabelian Baumslag--Solitar groups \cite{Mandel-Ushakov:2023b}, and
the lamplighter group \cite{Ushakov-Weiers:2023}.
Existence of a finitely generated abelian-by-cyclic group with undecidable spherical equations was established in \cite{Dong:2025}.
It was claimed in \cite{Kharlampovich-Lopez-Miasnikov:2020} 
that the Diophantine problem for systems of equations over $\MZ_2 \wr \MZ$ 
is decidable. However, a gap was later discovered in the proof, 
and the decidability of both systems and single equations 
over $\MZ_2 \wr \MZ$ remains an open question.
A brief survey on the solvability of equations in wreath products of groups 
can be found in \cite{Bartholdi-Dong-Pernak-Waechter:2024}.
Recently, O.~Kharlampovich and A.~Miasnikov 
\cite{KharlampovichMiasnikov2025}
showed that
the Diophantine problem in iterated wreath products 
of free abelian groups is undecidable.

Applications of quadratic equations (particularly spherical equations) 
to cryptography have also been explored. In particular, 
\cite{Ushakov:2024} establishes connections
between computational group theory and lattice-based cryptography 
via spherical equations over a certain class of finite metabelian groups.

\subsection{Our contribution}

As mentioned above, it remains an open question whether the Diophantine problem
for arbitrary single equations over the lamplighter group is decidable.
Nevertheless, decidability has been established for certain restricted classes
of equations, such as quadratic ones, \cite{Ushakov-Weiers:2023}.
In this paper, we show that one-variable equations constitute another such class:
they are decidable, with worst-case time complexity bounded by an
exponential function.
Moreover, for a ``typical'' one-variable equation, the existence of a solution
can be decided in nearly quadratic time.

\subsection{Study of one-variable equations}

One-variable equations can be solved efficiently in finitely generated abelian groups.  

R.~Lyndon \cite{Lyndon:1960(2)} was the first to study one-variable 
equations over free groups,
characterizing their solution sets in terms of parametric words.  
These parametric words were later simplified by A. Lorents 
\cite{Lorents:1963,Lorents:1968} and K.~Appel \cite{Appel:1968}.  
However, Lorents announced his results without proof, and Appel's published proof contains a gap 
(see \cite{Chiswell-Remeslennikov:2000}).  
A complete proof was subsequently provided by I.~Chiswell and V.~Remeslennikov.  
Solution sets of one-variable equations over free groups
were also studied via context-free languages 
\cite{Gilman-Myasnikov:2004}, leading to a high-degree polynomial-time algorithm 
\cite{Bormotov-Gilman-Miasnikov:2008}.  
More recently, R.~Ferens and A.~Je\.{z} \cite{Ferens-Jez:2021} 
described a cubic-time algorithm for the same problem.  

One-variable equations over nilpotent groups were studied by N.~Repin \cite{Repin:1984}.  
He proved that for any \(c \ge 10^{20}\), there is no algorithm to recognize the solvability of equations 
with one indeterminate in free nilpotent groups of class \(c\).  
Nevertheless, such an algorithm exists for finitely generated nilpotent groups of class $2$ 
(see \cite{Levine2022_virtually_class2_nilpotent,DuchinLiangShapiro2015_equations_nilpotent_groups}).  
Repin also constructed a finitely generated nilpotent group of class $3$ 
for which no algorithm exists to recognize the solvability of equations in one unknown.

\subsection{Model of computation and data representation}
\label{se:model-computations}

Below we show that the problem of deciding a one-variable equation
over the lamplighter group reduces to (multiple instances of)
divisibility of polynomials over the finite field $\MZ_2$.
Accordingly, our complexity estimates reflect the cost of polynomial
divisibility over $\MZ_2$.
Since all arithmetic operations take place in a field of size two,
individual field operations have constant cost.
We therefore work in the standard RAM model of computation, with
polynomials represented as described in Section~\ref{se:representation}.

We use the standard Bachmann--Landau asymptotic notation
(Big-O and related notation).
In particular, for functions
\(f,g:\mathbb N\to\mathbb R_{\ge 0}\),
the relation $f(n)$ is $O(g(n))$
means that there exist constants \(C>0\) and \(n_0\in\mathbb N\) such that
\[
f(n)\le C\,g(n)
\]
for all \(n\ge n_0\).

We also use the soft-O notation \(\wO\).
More precisely, for functions
\(f,g:\mathbb N\to\mathbb R_{\ge 0}\), we write
$f(n)$ is $\widetilde{O}(g(n))$
if there exists a constant \(k\geq 0\) such that
$f(n)$ is $O\bigl(g(n)(\log n)^k\bigr)$.
Equivalently, \(\widetilde{O}\) suppresses polylogarithmic factors 
in the input size \(n\).

\subsection{Outline}

We briefly describe the structure of the paper and the role of each section.

In Section~\ref{se:Laurent}, we review basic notions related to Laurent
polynomials over $\MZ_2$, with particular emphasis on their data representation
and the computational complexity of the fundamental arithmetic operations.
This section establishes the computational model used throughout the paper
and provides complexity bounds that are later invoked repeatedly.

In Section~\ref{se:lamplighter}, we recall the definition of the lamplighter
group $L_2=\MZ_2\wr\MZ$ and fix notation for its standard representation.
We also introduce one-variable equations over $L_2$ and reformulate them in
a form amenable to algorithmic analysis.

In Section~\ref{se:equations-to-divisibility}, we introduce the notion
of a \emph{$\delta$-parametric polynomial}, namely an expression of the form
\[
\sum_{i=s}^t f_i(z) z^{i\delta},
\]
where $f_s(z),\dots,f_t(z)\in\MZ_2[z^\pm]$ are fixed Laurent polynomials
and $\delta$ is a parameter.
For a given one-variable equation $w=1$ over $L_2$, we construct
two $\delta$-parametric polynomials
$\blf=\den(w)$ and $\blg=\num(w)$,
and show that decidability of $w=1$ reduces to determining whether
$\blg_\delta$ is divisible by $\blf_\delta$
for some $\delta\in\MZ$.
This reduction naturally splits the problem into two cases according to
whether $\sigma_x(w)$ vanishes, and identifies divisibility as the main
computational task.

In Section~\ref{se:num-den},
we investigate structural properties of the $\delta$-parametric polynomials
$\blf=\den(w)$ and $\blg=\num(w)$.
We show that these polynomials admit a geometric interpretation in terms of
finite subsets of the $\MZ^2$ grid, and develop a word-tracing procedure
that allows them to be computed efficiently.
We also establish several bounds on their size and support that will be used
throughout the remainder of the paper.

In Section~\ref{se:divisibility-decidable}, we study the divisibility
problem for $\delta$-parametric polynomials and establish its decidability.
The main idea is to replace the parametric monomial $z^\delta$
in $\blf$ and $\blg$ by an auxiliary variable $y$, thereby obtaining
bivariate polynomials $F(z,y)$ and $G(z,y)$.
We then perform symbolic division of $G$ by $F$ over the field
$\MZ_2(z)$ and analyze the resulting quotient and remainder.
This naturally leads to two cases, depending on whether 
the remainder vanishes.
In each case, we derive explicit bounds on the parameter $\delta$,
which yield decision procedures for divisibility and ultimately show
that the problem belongs to $\NP$.

Finally, in Section \ref{se:generic-complexity} 
we analyze the complexity of the resulting algorithm and show that,
although the worst-case running time is exponential, the problem admits
a nearly quadratic-time solution in the generic case.
This establishes decidability of one-variable equations over the lamplighter
group and demonstrates a strong contrast between worst-case and typical
behavior.

\section{Preliminaries: Laurent polynomials}
\label{se:Laurent}

Let $R$ be a commutative ring with unity.
A \emph{Laurent polynomial} over $R$ is an expression of the form
$$
\sum_{i=-\infty}^\infty a_i z^i\ \ \ (a_i\in R),
$$
where all but finitely many \emph{coefficients} $a_i$ 
are trivial. The sum of two polynomials is defined by
$$
\sum_{i=-\infty}^\infty a_i z^i +
\sum_{i=-\infty}^\infty b_i z^i = 
\sum_{i=-\infty}^\infty (a_i + b_i) z^i
$$
and the product of two polynomials is defined by
$$
\rb{\sum_{i=-\infty}^\infty a_i z^i} \cdot
\rb{\sum_{i=-\infty}^\infty b_i z^i} = 
\sum_{i=-\infty}^\infty c_i z^i,
\mbox{ where } c_i=\sum_{j=-\infty}^\infty a_j b_{i-j}.
$$
The set of all Laurent polynomials, denoted $R[z^\pm]$,
equipped with $+$ and $\cdot$ defined above is a ring, see \cite{CoxLittleOSheaIVA}.
Throughout the paper the ring of coefficients $R$ is
the finite field $\MZ_2=\{0,1\}$ of order $2$.

We say that a Laurent polynomial is \emph{trivial} if all its coefficients are zero.
For a nontrivial $f(z) = \sum_{i=-\infty}^\infty a_i z^i$ define
\begin{itemize}
\item 
$\deg(f) = \max \Set{i}{a_i\ne 0}$ called the \emph{degree} of $f$,
\item 
$\ord(f) = \min \Set{i}{a_i\ne 0}$ called the \emph{order} of $f$.
\end{itemize}
The monomial $a_{\deg(f)} z^{\deg(f)}$ in $f$
is called the \emph{leading monomial}
and $a_{\ord(f)} z^{\ord(f)}$ is called the \emph{trailing monomial}.
For the trivial polynomial $f=0$,
$\ord(f)$ and $\deg(f)$ are not defined
and we denote this by writing $\deg(f)=\varnothing$ and $\ord(f)=\varnothing$.
Define 
$$
\size(f) =
\begin{cases}
\max(|\deg(f)|,|\ord(f)|) & \mbox{if } f\ne 0,\\
0 & \mbox{if } f=0.\\
\end{cases}
$$

\subsection{Division in $\MZ_2[z^\pm]$}

To divide $f(z)$ by $g(z)\ne 0$ means to find $q(z),r(z) \in\MZ_2[z^\pm]$
satisfying $f(z) = g(z)\cdot q(z) + r(z)$, where
\begin{equation}\label{eq:r}
r=0\ \ \mbox{ or }\ \ 
\ord(f)\le \ord(r)\le \deg(r) <
\ord(f)+[\deg(g)-\ord(g)].
\end{equation}
If $r=0$, then we say that $g$ \emph{divides} $f$ and write $g \mid f$.
We call $q(z)$ and $r(z)$ the \emph{quotient} and \emph{remainder}, respectively.
Denote by $g\ \%\ f$ the remainder of division of $f$ by $g$.
There is a close relationship between division in $\MZ_2[z^{\pm}]$
and division in $\MZ_2[z]$.  
Indeed, for any $f,g,r \in \MZ_2[z^{\pm}]$ with $g \neq 0$,
we have
\begin{align*}
r = f\ \%\ g\ \  \Leftrightarrow\ \ & 
f=gq+r \mbox{ for some }q\in \MZ_2[z^\pm] \mbox{ satisfying \eqref{eq:r}}\\
\ \  \Leftrightarrow\ \ & 
\underbrace{f\cdot z^{-\ord(f)}}_{F\in\MZ_2[z]} = 
\underbrace{g \cdot z^{-\ord(g)}}_{G\in\MZ_2[z]}
\cdot 
\underbrace{q \cdot z^{-\ord(q)}}_{Q\in\MZ_2[z]} + 
\underbrace{r\cdot z^{-\ord(f)}}_{R}\\
& \scalebox{0.8}{(satisfying $R=0$ or $0\le \ord(R),\ \ \deg(R)<\deg(G)$)} \\
\ \  \Leftrightarrow\ \ & 
R = F\ \%\ G \mbox{ in } \MZ_2[z].
\end{align*}
In particular, the following holds.

\begin{lemma}\label{Lemma:equivalence of Laurent division}
$g(z)$ divides $f(z)$ in $\MZ_2[z^\pm]$ $\ \ \Leftrightarrow\ \ $
$G(z)$ divides $F(z)$ in $\MZ_2[z]$.
\end{lemma}

\subsection{Representation for Laurent polynomials}
\label{se:representation}

Recall that in \cite{GathenGerhard2003} a polynomial over a ring $R$
is represented by an array whose $i$th entry is the coefficient $a_i$.
In a similar way, we represent a nontrivial Laurent polynomial
\[
f = \sum_{i=-\infty}^{\infty} a_i z^i
\]
with coefficients in $\MZ_2$ by the pair $(\deg(f), w_f)$, where $\deg(f)$ is given in unary and
$w_f \in \{0,1\}^\ast$ is a bit string defined by
\[
w_f =
\begin{cases}
a_{\ord(f)} \dots a_{\deg(f)}, & \text{if } f \neq 0, \\
\varepsilon, & \text{if } f = 0.
\end{cases}
\]
The zero polynomial $f = 0$ is represented by the pair $(0, w_0)$.
Observe that the value $\size(f)$, defined as above, accurately reflects
the size of the representation of $f$ -- the number of bits required to store
the pair $(\deg(f), w_f)$ is $\Theta(\size(f))$.
Note also that, by definition, if $w_f$ is nonempty, then it both starts and
ends with $1$.
We refer to such a representation as \emph{strict}.

In some situations it is convenient to use a \emph{non-strict} representation
of a polynomial $f$ defined by a pair $(D, w)$, where $\deg(f) \le D$,
$w = a_d \dots a_D$, and $d \le \ord(f)$.
Unlike the strict representation, a non-strict representation allows
the word $w$ to begin and end with $0$.
We write $\pres(f) = (D, w)$ if the pair $(D, w)$ represents the polynomial $f$.

\subsection{Complexity of operations in $\MZ_2[z^\pm]$}
\label{se:complexity-operations}

Let $R$ be any commutative ring with unity.
Addition, subtraction, multiplication, and division
for polynomials from $R[z]$ can be performed in linear or nearly linear time.
In more detail, the following holds.

\begin{theorem}[{{\cite[Section 2.2]{GathenGerhard2003}}}]
\label{th:plus}
Polynomials of degree less than $n$  can be added using at most $O(n)$
arithmetic operations in $R$.   
\end{theorem}

\begin{theorem}[{{\cite[Theorem 8.23]{GathenGerhard2003}}}]
\label{th:times}
Polynomials of degree less than $n$ can be multiplied using at most 
$(18+72 \log_3 2)n\log n \log\log n + O(n\log n)$
or $63.43n\log n \log\log n + O(n\log n)$
arithmetic operations in $R$.
\end{theorem}

Denote by $M(n)$ the number of arithmetic operations in $R$
required to multiply two arbitrary polynomials in $R[x]$ 
of degree less than $n$.
By Theorem \ref{th:times}, $M(n)\in O(n\log n \log \log n )$.

\begin{theorem}[{{\cite[Section 9.1]{GathenGerhard2003}}}]
\label{th:div}
Division with remainder of a polynomial $f\in R[x]$
by a monic polynomial $g\in R[x]$
of degree less than $n$ can be done using 
$O(M(n))$ arithmetic operations in $R$.
\end{theorem}

Addition, subtraction, multiplication, and division of Laurent polynomials 
in $\MZ_2[z^{\pm}]$ can be reduced to the corresponding operations 
on ordinary polynomials in $\MZ_2[z]$.
Indeed, consider $f(z),g(z)\in \MZ_2[z^{\pm}]$. Define 
$$
o_f=\ord(f),\ o_g=\ord(g),\ \mbox{ and } m=-\min(o_f,o_g),
$$
and notice that in this notation
\begin{align}
\label{eq:plus}
f(z) \pm g(z) &= (f(z)z^{-m} \pm g(z)z^{-m})z^m\\
\label{eq:times}
f(z) \cdot g(z) &= (f(z)z^{-o_f} \cdot g(z)z^{-o_g})z^{o_f+o_g}\\
\label{eq:div}
g(z) \% f(z) &= (g(z)z^{-o_g}\ \%\ f(z)z^{-o_f})z^{o_f}.
\end{align}

\begin{proposition}\label{pr:Laurent-complexity}
Let $f(z),g(z)\in \MZ_2[z^{\pm}]$ and 
$
n = \max\{\size(f),\size(g)\}
$.
\begin{itemize}
\item[(a)]
There is an algorithm that computes $f(z)\pm g(z)$ in $O(n)$ time.
\item[(b)]
There is an algorithm that computes $f(z)\cdot g(z)$ in $O(n\log n\log \log n)$ time.
\item[(c)]
There is an algorithm that computes $g(z)\ \%\ f(g)$ in $O(n\log n\log \log n)$ time.
\item[(d)]
There is an algorithm that decides 
whether $g(z)\ne 0$ divides $f(z)$ in $O(n\log n\log \log n)$ time.
\end{itemize}
\end{proposition}

\begin{proof}
Theorems~\ref{th:plus}, \ref{th:times}, and \ref{th:div} allow us to compute
\[
f(z) z^{-m} \pm g(z) z^{-m}, \qquad
f(z) z^{-o_f} \cdot g(z) z^{-o_g}, \qquad \text{and} \qquad
g(z) z^{-o_g} \bmod f(z) z^{-o_f},
\]
in time $O(n)$, $O(n \log n \log \log n)$, and
$O(n \log n \log \log n)$, respectively.
Multiplying the resulting polynomials by $z^{m}$, $z^{o_f + o_g}$, and $z^{o_f}$,
respectively, yields \eqref{eq:plus}, \eqref{eq:times}, and \eqref{eq:div}.
\end{proof}

In Section \ref{se:divisibility-decidable}, we also employ 
the \emph{synthetic division} algorithm (see \cite[Section 2.5]{GathenGerhard2003}), 
whose complexity can be estimated as $O(n^2)$ arithmetic operations in $R$.




\section{Preliminaries: the lamplighter group}
\label{se:lamplighter}

Define the support of a function $f\colon \MZ\to \MZ_2$ as a set
$$
\supp(f) =\Set{i\in \MZ}{f(i)\ne 0}.
$$
Define a set 
$$
\MZ_2^\MZ  = \Set{f\colon \MZ\to \MZ_2}{|\supp(f)|<\infty}
$$
and a binary operation $+$ on $\MZ_2^\MZ$ which for
$f,g\in \MZ_2^\MZ$ produces $f+g\in \MZ_2^\MZ$ defined by
$$
(f+g)(z) = f(z)+g(z) \ \mbox{ for } \ z\in \MZ.
$$
For a nontrivial $f$ it will be convenient to define
$$
m(f)=\min_{z\in \supp(f)} z
\ \ \mbox{ and }\ \ 
M(f)=\max_{z\in \supp(f)} z.
$$
For $f\in \MZ_2^\MZ$ and $b\in \MZ$, define $f^b\in \MZ_2^\MZ$ by
$$
f^b(z)=f(z+b) \
\mbox{ for } \ z\in \MZ,
$$
which is a right $\MZ$-action on $\MZ_2^\MZ$ because $f^0=f$ and 
$f^{(b_1+b_2)}(z)=f(z+b_1+b_2)=(f^{b_1})^{b_2}(z)$.
Hence we can consider a semi-direct product $\MZ \ltimes \MZ_2^\MZ$
equipped with the operation
\begin{equation}\label{eq:semi-product}
(\delta_1,f_1)(\delta_2,f_2)=
(\delta_1+\delta_2,f_1^{\delta_2}+f_2).
\end{equation}
The group $\MZ \ltimes \MZ_2^\MZ$ is called the
\emph{restricted wreath product} of $\MZ_2$ and $\MZ$ and is denoted by
$\MZ_2\wr\MZ$. The group $\MZ_2\wr\MZ$ is also well-known as the \textit{lamplighter group}(\cite{Ushakov-Weiers:2023}), 
as it can be viewed as an infinite set of lamps (each lamp indexed by an element 
of $\MZ$), with each lamp either on or off, and a lamplighter positioned at some lamp.
Given some element $(\delta,f)\in \MZ_2\wr\MZ$, $f\in\MZ_2^\MZ$ represents the 
configuration of illuminated lamps and $\delta\in\MZ$ represents the position 
of the lamplighter. Henceforth, we denote $\MZ_2\wr\MZ$ as $L_2$.

There is a natural (abelian group) isomorphism between $\MZ_2^\MZ$ and the 
ring  $\MZ_2[z^\pm]$ of Laurent polynomials
with coefficients in $\MZ_2$ that maps 
$f\in \MZ_2^\MZ$ to $\sum_{i=-\infty}^\infty f(i)z^i$.
Hence, the group $L_2$ can be viewed 
as $\MZ\ltimes \MZ_2[z^\pm]$ with the $\MZ$-action on $\MZ_2[z^\pm]$
defined by
$$
f^\delta = f\cdot z^{-\delta}.
$$

It is easy to see that $L_2$ is a two-generated group.
In particular, it can be generated by the elements
$a=(0,\mathbf{1}_0)$ and $t = (1,0)$, where $\one_0 \in \MZ_2^\MZ$ is defined by
$$
\one_0(x)=
\begin{cases}
1& \mbox{ if } x=0\\
0& \mbox{ if } x\ne0\\
\end{cases}
$$
and $0$ in the second component of $t$ is the function identical to $0$.
G.~Baumslag in \cite{Baumslag:1961} proved that $L_2$ is not finitely presented,
but can be defined using the following infinite presentation:
$$
L_2 \simeq \gpr{a,t}{a^2, [a^{t^i}, a^{t^j}](i,j \in \MZ)}.
$$

\begin{lemma}\label{le:word-to-pair}
It takes $O(|w|)$ time to compute the pair $(\delta,f) \in \MZ\ltimes \MZ_2[z^\pm]$
for a given group-word $w=w(a,t)$.
\end{lemma}

\begin{proof}
One can process $w$ letter-by-letter from the right to the left using the formula
\eqref{eq:semi-product} with $(\delta_1,f_1) = (0,\mathbf{1}_0)^\pm$ or $(1,0)^\pm$
that can be handled using formulae
$$
(0,\mathbf{1}_0)^\pm \cdot (\delta_2,f_2) = (\delta_2,f_2+z^{\delta_2})
\ \ \mbox{ and }\ \ 
(1,0)^\pm \cdot (\delta_2,f_2) = (\delta_2\pm 1,f_2).
$$
This can be done in $O(|w|)$ time.
\end{proof}

Finally, we formulate the Diophantine problem for one-variable equations over $L_2$.

\begin{algproblem}
\problemtitle{\textsc{Diophantine problem for one-variable equations over $L_2$}, $(\DP_1(L_2))$.}
\probleminput{A group word $w=w(a,t,x)\in F(a,t,x)$.}
\problemquestion{Is there $s=s(a,t)\in F(a,t)$ satisfying $w(a,t,s(a,t))=1$ 
in $L_2$?}
\end{algproblem}

For a word $w=w(a,t,x)\in F(a,t,x)$ define
\begin{itemize}
\item 
$\sigma_a(w) = $ the sum of exponents for $a$ in $w$,
\item 
$\sigma_t(w) = $ the sum of exponents for $t$ in $w$,
\item 
$\sigma_x(w) = $ the sum of exponents for $x$ in $w$.
\end{itemize}
$\sigma_a$, $\sigma_t$, and $\sigma_x$
are computable in linear time $O(|w|)$ for a given $w=w(a,t,x)$.

\subsection{Useful formulae}

For any $(\delta,f),(\delta_1,f_1) \in \MZ_2\wr \MZ$
the following holds:
\begin{align*}
(\delta,f)^{-1} &= (-\delta,-f^{-\delta})
&& \mbox{(inverse)}\\
(\delta,f)^{-1} (\delta_1,f_1) (\delta,f) &= 
(\delta_1,(1-z^{-\delta_1})f+z^{-\delta}f_1)
&& \mbox{(conjugation)}\\
[(\delta, f), (\delta_1, f_1)] &= (0, -f - f_1^ {\delta} + f^{\delta_1} + f_1)
&& \mbox{(commutator)}\\
&= (0, -f(1-z^{-\delta_1}) + f_1(1-z^{-\delta}))\\
&= (0, f(1-z^{-\delta_1}) + f_1(1-z^{-\delta})).
\end{align*}

\section{One-variable equations: translation to parametric polynomials}
\label{se:equations-to-divisibility}

In this section we introduce \emph{$\delta$-parametric polynomials}
and prove that the Diophantine problem $\DP_1(L_2)$
for one-variable equations over $L_2$ 
can be reduced,
in polynomial-time,
to divisibility problem for $\delta$-parametric polynomials.

\begin{lemma}
It requires linear time to translate a word $w(a,t,x)$ into the form 
\begin{equation}\label{eq:one-var-eq}
c_0 x^{\varepsilon_1} c_1
\dots 
c_{k-1} x^{\varepsilon_k} c_k = 1,
\end{equation}
where $c_i=(\delta_i,f_i)\in \MZ \ltimes \MZ_2^\MZ$ are constants
and $\varepsilon_i=\pm 1$.
\end{lemma}

\begin{proof}
It is straightforward to rewrite $w(a,t,x)$ as an alternating product
$$
w_0(a,t) x^{\varepsilon_1} w_1(a,t)
\dots 
w_{k-1}(a,t) x^{\varepsilon_k} w_k(a,t),
$$
and then apply Lemma~\ref{le:word-to-pair} to each word $w_i(a,t)$ individually.
\end{proof}

Hence, we may assume that the equation is given in the form \eqref{eq:one-var-eq}.
For $i=0,\dots,k$ define
\begin{align*}
\arraycolsep=3pt
\begin{array}{rcl}
\alpha_i
&=&\delta_{i+1}+\dots+\delta_{k}\\
\beta_i
&=&\varepsilon_{i+1}+\dots+\varepsilon_{k}
\end{array}
&&
\gamma_i =
\begin{cases}
\beta_i &\mbox{if } \varepsilon_i=1,\\
\beta_i-1 &\mbox{if } \varepsilon_i=-1.\\
\end{cases}
\end{align*}
With this notation, substituting $x=(\delta,f)$ into the left-hand 
side of \eqref{eq:one-var-eq} yields
\begin{align}\nonumber
(\delta_0,f_0)
(\delta,f)^{\varepsilon_1}
\dots
(\delta,f)^{\varepsilon_k}
(\delta_k,f_k) 
&= 
\Bigl(\sum_{i=0}^k \delta_i + \delta\sum_{i=1}^k \varepsilon_i,
\ \ 
\sum_{i=0}^k f_i^{\alpha_i+\beta_i\delta}+
\sum_{i=1}^k \varepsilon_i f^{\alpha_{i-1}+\gamma_i\delta}
\Bigr)\\
\nonumber
&= 
\Bigl(\sum_{i=0}^k \delta_i + \delta\sum_{i=1}^k \varepsilon_i,
\ \ 
\sum_{i=0}^k f_i z^{-\alpha_i-\beta_i\delta}+
f
\sum_{i=1}^k \varepsilon_i z^{-\alpha_{i-1}-\gamma_i\delta}
\Bigr)\\
\label{eq:rhs}
&= 
\Bigl(\underbrace{\sum_{i=0}^k \delta_i}_{t_w}
+ \delta\underbrace{\sum_{i=1}^k \varepsilon_i}_{x_w},
\ \ 
\underbrace{\sum_{i=0}^k f_i z^{-\alpha_i-\beta_i\delta}}_{\num(w)}+
f
\underbrace{\sum_{i=1}^k z^{-\alpha_{i-1}-\gamma_i\delta}}_{\den(w)}
\Bigr).
\end{align}
In the last expression, the coefficients $\varepsilon_i = \pm 1$ are omitted
from the definition of $\den(w)$, since $-1 = 1$ in $\MZ_2$.
We denote by $t_w$, $x_w$, $\num(w)$, and $\den(w)$ the four components
of~\eqref{eq:rhs}, namely
\begin{align*}
t_w&=\sum_{i=0}^k \delta_i, &
x_w&=\sum_{i=0}^k \varepsilon_i, &
\num(w)&=\sum_{i=0}^k f_iz^{-\alpha_i-\beta_i\delta},&
\den(w)&=\sum_{i=1}^k z^{-\alpha_{i-1}-\gamma_i\delta}.
\end{align*}

\begin{lemma}
$\sigma_x(w) = x_w$ and $\sigma_t(w) = t_w$ for every $w=w(a,t,x)$.
\end{lemma}

\begin{proof}
The sums $\sum_{i=1}^k \varepsilon_i$ and $\sum_{i=0}^k \delta_i$ 
record the total exponents of $x$ and $t$
in the word $w(a,t,x)$. Thus, they compute
$\sigma_x(w)$ and $\sigma_t(w)$ respectively.
\end{proof}

\begin{lemma}\label{le:C1-C2}
$x=(\delta,f)$ satisfies  \eqref{eq:one-var-eq} 
$\ \Leftrightarrow\ $
the following holds for $\delta$ and $f$:
\begin{itemize}
\item[(C1)]
$\sum_{i=0}^k \delta_i + \delta\sum_{i=1}^k \varepsilon_i = 0$,
\item[(C2)]
$\sum_{i=0}^k f_i z^{-\alpha_i-\beta_i\delta}+
f\sum_{i=1}^k z^{-\alpha_{i-1}-\gamma_i\delta} = 0$.
\end{itemize}
\end{lemma}

\begin{proof}
Indeed, $x=(\delta,f)$ satisfies \eqref{eq:one-var-eq} if and only if
\eqref{eq:rhs} evaluates to $(0,0)$,
which yields the conditions (C1) and (C2).
\end{proof}

\subsection{Parametric polynomials}

The expressions $\num(w)$ and $\den(w)$ are Laurent polynomials
in which 
the parameter $\delta$ is not specified.
We refer to such expressions as \emph{$\delta$-parametric polynomials}.
Formally, a \emph{$\delta$-parametric polynomial} is an expression of the form
\begin{equation}\label{eq:param-poly}
\sum_{i=s}^t f_i(z) z^{i\delta},
\end{equation}
where $f_s(z),\dots,f_t(z) \in \MZ_2[z^\pm]$ are fixed Laurent polynomials
and $\delta$ is a \emph{parameter}, that is, a variable ranging over $\MZ$.
We use boldface letters, such as $\blf$ or $\blg$, 
to denote
$\delta$-parametric polynomials.
To \emph{instantiate} a $\delta$-parametric polynomial means to
substitute the parameter $\delta$ by a specific integer value.
The result of instantiating $\blf$ at $\delta \in \MZ$
is a Laurent polynomial, denoted by $\blf_\delta$.

For a $\delta$-parametric polynomial $\blf = \sum_{i=s}^t f_i(z) z^{i\delta}$
(where $f_s\ne0$, $f_t\ne0$)
define its \emph{$\delta$-degree} and \emph{$\delta$-order} as
$$
\deg_\delta(\blf) = t\ \ \mbox{ and }\ \ \ord_\delta(\blf) = s.
$$


\subsection{Criteria for solvability}

Condition~(C1) gives rise to two cases for the equation $w(a,t,x) = 1$,
namely $x_w \neq 0$ and $x_w = 0$.

\begin{proposition}
\label{pr:xw-nonzero}
Suppose $x_w \neq 0$. Then the equation $w(a,t,x) = 1$ has a solution
if and only if $x_w \mid t_w$ and, for
\begin{equation}
\label{eq:delta}
\delta = -\frac{t_w}{x_w},
\end{equation}
one of the following conditions holds:
\begin{itemize}
\item
$\den_\delta(w) \neq 0$ and $\den_\delta(w) \mid \num_\delta(w)$, or
\item
$\den_\delta(w) = 0$ and $\num_\delta(w) = 0$.
\end{itemize}
\end{proposition}

\begin{proof}
When $x_w \neq 0$, condition~(C1) uniquely determines the value of $\delta$
as in~\eqref{eq:delta}. Moreover, condition~(C2) then uniquely determines
the polynomial $f$:
\begin{equation}
\label{eq:f}
f = -\frac{\num_\delta(w)}{\den_\delta(w)} = \frac{\num_\delta(w)}{\den_\delta(w)}
\quad \text{if } \den_\delta(w) \neq 0,
\end{equation}
or $f$ can be any polynomial if $\den_\delta(w) = 0$ and $\num_\delta(w) = 0$.
These yield the stated conditions.
\end{proof}


In Proposition~\ref{Prop:algorithm-equations} we prove that
there is a polynomial time algorithm deciding
whether a given equation $w(a,t,x)=1$ satisfying $x_w \neq 0$
has a solution.

When $x_w=\sum_{i=1}^k \varepsilon_i = 0$, condition~(C1), when consistent,
does not determine $\delta$, which may then range over all integers.
This makes the problem significantly more difficult.

\begin{proposition}\label{pr:xw-zero}
Suppose $x_w = 0$. Then the equation $w(a,t,x) = 1$ has a solution
if and only if $t_w = 0$ and there exists $\delta \in \MZ$ such that
one of the following conditions holds:
\begin{itemize}
\item[(C2')]
$\den_\delta(w) \neq 0$ and $\den_\delta(w) \mid \num_\delta(w)$, or
\item[(C2'')]
$\den_\delta(w) = 0$ and $\num_\delta(w) = 0$.
\end{itemize}
\end{proposition}

\begin{proof}
If $x_w = 0$, then condition~(C1) is satisfied precisely when $t_w = 0$.
Moreover, the disjunction of (C2') and (C2'')
is equivalent to condition~(C2) for some value of $\delta$.
\end{proof}

\subsection{Divisibility problems for parametric polynomials}
\label{se:DIV-problem}

Here we define two versions of the divisibility problem for parametric polynomials.

\begin{algproblem}
 \problemtitle{\textsc{Divisibility for $\delta$-parametric polynomials} $(\DIV(\blf,\blg))$.}
  \probleminput{$\delta$-parametric polynomials $\blf$ and $\blg$.}
  \problemquestion{Does there exist $\delta \in \MZ$ such that $\blf_\delta \mid \blg_\delta$?}
\end{algproblem}

An integer $\delta\in\MZ$ satisfying $\blf_\delta\mid \blg_\delta$ is called
a \emph{witness} for the instance $(\blf,\blg)$ of $\DIV$.

\begin{proposition}\label{prop:main-equivalence}
The problem of deciding whether a given one-variable equation
$w(a,t,x) = 1$ satisfying $\sum_{i=1}^k \varepsilon_i = 0$
admits a solution is polynomial-time reducible to determining whether
the pair $(\den(w), \num(w))$ is a positive instance of $\DIV$.
\end{proposition}

\begin{proof}
First, compute $\sigma_t(w) = t_w$ in $O(|w|)$ time and check whether $t_w = 0$.

By Lemma~\ref{le:trivial-den-num}, the existence of 
$\delta \in \MZ$ such that $\den_\delta(w) = 0$ and $\num_\delta(w) = 0$
(condition (C2'') in Proposition~\ref{pr:xw-zero}) is equivalent to one of the following:
\begin{itemize}
\item $\den(w) = 0$ and $\num(w) = 0$, or
\item $\den_\delta(w) = 0$ and $\num_\delta(w) = 0$ for some
$\delta$ with $- |w| \le \delta \le |w|$.
\end{itemize}
By Proposition~\ref{pr:compute-num-den}, if $|\delta| \le |w|$, then
$\num_\delta(w)$ and $\den_\delta(w)$ can be computed in $O(|w|^2)$ time.
Hence, both conditions can be checked in polynomial time for any given word $w$.

Therefore, when $x_w = 0$, deciding whether the equation $w(a,t,x) = 1$
admits a solution reduces to determining whether there exists
$\delta \in \MZ$ such that $\den_{\delta}(w)$ divides $\num_{\delta}(w)$
(condition (C2') of Proposition~\ref{pr:xw-zero}).
This corresponds to checking whether the pair 
$(\den(w), \num(w))$ constitutes a positive instance of the divisibility problem.
\end{proof}

The problem $\DIV$ naturally splits into two cases, according to the sign 
of $\delta$: $\delta > 0$ and $\delta < 0$.
We claim that both cases can be reduced to instances of the following problem.

\begin{algproblem}
 \problemtitle{\textsc{Positive divisibility for $\delta$-parametric polynomials} $(\DIV_+(\blf,\blg))$.}
  \probleminput{$\delta$-parametric polynomials $\blf$ and $\blg$.}
  \problemquestion{Does there exist $\delta > 0$ such that $\blf_\delta \mid \blg_\delta$?}
\end{algproblem}

An integer $\delta > 0$ satisfying $\blf_\delta \mid \blg_\delta$
is called a \emph{witness} for the instance $(\blf, \blg)$ of $\DIV_+$.
For a $\delta$-parametric polynomial $\blf=\sum_{i=s}^t f_i(z) z^{i\delta}$
we define the \emph{reverse} of $\blf$ as
$$
\rev(\blf) = \sum_{i=s}^t f_i(z) z^{-i\delta}
$$
obtained by replacing each occurrence of $\delta$ with $-\delta$.
Clearly, $[\rev(\blf)]_\delta = \blf_{-\delta}$ for every $\delta\in\MZ$.
Therefore,
\begin{align*}
\exists \delta < 0 \mbox{ such that } 
\blf_\delta \mid \blg_\delta
&\ \ \Leftrightarrow\ \ 
\exists \delta < 0 \mbox{ such that } 
[\rev(\blf)]_{-\delta} \mid [\rev(\blg)]_{-\delta}\\
&\ \ \Leftrightarrow\ \ 
\exists \delta > 0 \mbox{ such that } 
[\rev(\blf)]_{\delta} \mid [\rev(\blg)]_{\delta}.
\end{align*}
It will be convenient to assume that 
$\DIV(\blf,\blg))$ and $\DIV_+(\blf,\blg))$ define a true/false value
indicating whether $(\blf, \blg)$
is a positive/negative instance of the problem.
As a corollary we obtain the following proposition.

\begin{proposition}\label{pr:div-div2}
For all parametric polynomials $\blf,\blg$ the following holds:
\begin{equation}\label{eq:div-div2}
\DIV(\blf,\blg)
=
\DIV_+(\blf,\blg)
\;\vee\;
(\blf_0\mid\blg_0)
\;\vee\;
\DIV_+(\rev(\blf),\rev(\blg)).
\end{equation}
\end{proposition}

\section{Properties of $\num_{\delta}$ and $\den_{\delta}$}
\label{se:num-den}

\subsection{Inductive formulae for $\num$ and $\den$}

Consider an equation $w(a,t,x)=1$ and the \emph{compact form} 
of the expression \eqref{eq:rhs} for the word $w$
\begin{equation}\label{eq:rhs-short}
(\delta_0,f_0)
(\delta,f)^{\varepsilon_1}
\dots
(\delta,f)^{\varepsilon_k}
(\delta_k,f_k) 
\ =\ 
(x_w\delta + t_w,\ \num(w) + f \den(w)).
\end{equation}
The polynomials $\den(w)$ and $\num(w)$ can be efficiently constructed 
by reading the word $w$ \textbf{from the right to the left}
as described by the following lemmas.

\begin{lemma}\label{le:p-num-den0}
For the empty word $\varepsilon$ we have
$x_\varepsilon=t_\varepsilon=0$ and 
$\num(\varepsilon) = \den(\varepsilon) =0$.
\end{lemma}

\begin{lemma}\label{le:p-num-den1}
If $w'=x\circ w$, then 
$t_{w'}=t_w$, $x_{w'}=x_w+1$,
$\num(w')=\num(w)$, and
$\den(w')=\den(w) + z^{-x_w \delta-t_w}$.
\end{lemma}

\begin{proof}
Multiplying the formula \eqref{eq:rhs-short} 
on the left by $x=(\delta,f)$ we get
$$
(\delta,f)\cdot(x_w\delta + t_w,\ \num(w) + f \den(w)) =
((\underbrace{x_w+1}_{x_{w'}})\delta + \underbrace{t_w}_{t_{w'}},\ 
\underbrace{\num(w)}_{\num(w')} + 
f(\underbrace{\den(w) +z^{-x_w\delta  -t_w}}_{\den(w')}))
$$
which gives the claimed formulae.
\end{proof}

\begin{lemma}\label{le:p-num-den2}
If $w'= x^{-1} \circ w$, then 
$t_{w'} = t_w$, $x_{w'} = x_w -1$,
$\num(w')=\num(w)$, and
$\den(w')=\den(w) + z^{-(x_w-1) \delta - t_w}$.
\end{lemma}

\begin{proof}
Multiplying the formula \eqref{eq:rhs-short} 
on the left by $x^{-1}=(-\delta,-f^{-\delta})$ we get
$$
(-\delta,-f^{-\delta})\cdot(x_w\delta + t_w,\ \num(w) + f \den(w)) =
((x_w-1)\delta + t_w,\ \num(w) + f(\den(w) +z^{-(x_w - 1)\delta -t_w}))
$$
which gives the claimed formulae.
\end{proof}

\begin{lemma}\label{le:p-num-den3}
If $w'=t\circ w$, then 
$t_{w'}=t_w+1$, 
$x_{w'}=x_w$,
$\num(w')=\num(w)$, and
$\den(w')=\den(w)$.
\end{lemma}

\begin{proof}
Multiplying the formula \eqref{eq:rhs-short} 
on the left by $t=(1,\blo)$ we get
$$
(1,\blo)\cdot(x_w\delta + t_w,\ \num(w) + f \den(w)) =
(x_w\delta + (t_w+1),\ \num(w) + f\den(w))
$$
which gives the claimed formulae.
\end{proof}

\begin{lemma}\label{le:p-num-den4}
If $w'= t^{-1} \circ w$, then 
$t_{w'}=t_w-1$, 
$x_{w'}=x_w$,
$\num(w')=\num(w)$, and
$\den(w')=\den(w)$.
\end{lemma}

\begin{proof}
Multiplying the formula \eqref{eq:rhs-short} 
on the left by $t^{-1}=(-1,\blo)$ we get
$$
(-1,\blo)\cdot(x_w\delta + t_w,\ \num(w) + f \den(w)) =
(x_w\delta + (t_w-1),\ \num(w) + f\den(w))
$$
which gives the claimed formulae.
\end{proof}

\begin{lemma}\label{le:p-num-den5}
If $w' = a \circ w$, then 
$t_{w'}=t_w$,
$x_{w'}=x_w$,
$\num(w')=\num(w)+ z^{-x_w \delta-t_w}$, and
$\den(w')=\den(w)$.
\end{lemma}

\begin{proof}
Multiplying the formula \eqref{eq:rhs-short} 
on the left by $a=(0,\one_0)$ we get
$$
(0,\one_0)\cdot (x_w\delta + t_w,\ \num(w) + f \den(w)) =
(x_w\delta + t_w,\ (\num(w) +z^{-x_w \delta-t_w}) + f \den(w)) 
$$
which gives the claimed formulae.
\end{proof}

\subsection{Magnus-type embedding of one-variable equations}

Consider the free abelian group $\MZ^2$, defined as
\[
\MZ^2 = \{ x^i t^j  \mid i,j \in \MZ \},
\]
equipped with the binary operation
\[
x^i t^j \cdot x^k t^l = x^{i+k} t^{j+l} .
\]
Let $\MZ_2[\MZ^2]$ denote the corresponding group ring, 
consisting of all finite linear combinations of elements $x^j t^i$ from $\MZ^2$ 
with coefficients in $\MZ_2$.
Every \emph{$\delta$-parametric polynomial} over $\MZ_2$ 
can be written as a finite sum of monomials
\[
\sum_{i=1}^n z^{b_i \delta + a_i}, \quad a_i, b_i \in \MZ.
\]
In fact, the set of all $\delta$-parametric polynomials forms a ring, and there is a natural ring isomorphism
to the group ring $\MZ_2[\MZ^2]$ given by
\[
z^{b \delta + a} \;\longmapsto\; x^b t^a.
\]
Hence, every $\delta$-parametric polynomial can be represented as an element of $\MZ_2[\MZ^2]$.
Let $M$ be the group of $3 \times 3$ lower-triangular matrices of the form
$$
\left[
\begin{array}{ccc}
p & 0 & 0 \\
\blf & 1 & 0 \\
\blg & 0 & 1 \\
\end{array}
\right]
$$
where $p \in \MZ^2$ and $\blf,\blg\in \MZ_2[\MZ^2]$.
Define a homomorphism $\alpha: F(a,t,x) \to M$ by 
\[
x \mapsto
X = 
\begin{bmatrix}
x^{-1} & 0 & 0 \\
0 & 1 & 0 \\
1 & 0 & 1 
\end{bmatrix},
\ \ \ \ 
t \mapsto
T =
\begin{bmatrix}
t^{-1} & 0 & 0 \\
0 & 1 & 0 \\
0 & 0 & 1 \\
\end{bmatrix},
\ \ \ \ 
a \mapsto
A = 
\begin{bmatrix}
1 & 0 & 0 \\
1 & 1 & 0 \\
0 & 0 & 1 
\end{bmatrix}
\]

\begin{proposition}\label{pr:num-den-matrices}
$\alpha(w)=
\left[
\begin{array}{ccc}
x^{-x_w}t^{-t_w} & 0 & 0 \\
\num(w) & 1 & 0 \\
\den(w) & 0 & 1 \\
\end{array}
\right]
$ for every $w\in F(a,t,x)$.
\end{proposition}

\begin{proof}
Induction on $|w|$.
If $|w| = 0$, then $\alpha(w)=I$ which agrees with Lemma \ref{le:p-num-den0}.
Then using the induction assumption, we deduce the following:
{\allowdisplaybreaks
\begin{align*}
\alpha(x\circ w)
&=
\begin{bmatrix}
x^{-1} & 0 & 0 \\
0 & 1 & 0 \\
1 & 0 & 1 
\end{bmatrix}
\cdot
\begin{bmatrix}
x^{-x_w}t^{-t_w} & 0 & 0 \\
\num(w) & 1 & 0 \\
\den(w) & 0 & 1 \\
\end{bmatrix}
=
\begin{bmatrix}
x^{-(x_w+1)}t^{-t_w} & 0 & 0 \\
\num(w) & 1 & 0 \\
\den(w)+x^{-x_w}t^{-t_w} & 0 & 1 \\
\end{bmatrix}\\
\alpha(t\circ w)
&=
\begin{bmatrix}
t^{-1} & 0 & 0 \\
0 & 1 & 0 \\
0 & 0 & 1 \\
\end{bmatrix}
\cdot
\begin{bmatrix}
x^{-x_w}t^{-t_w} & 0 & 0 \\
\num(w) & 1 & 0 \\
\den(w) & 0 & 1 \\
\end{bmatrix}
=
\begin{bmatrix}
x^{-x_w}t^{-(t_w+1)} & 0 & 0 \\
\num(w) & 1 & 0 \\
\den(w) & 0 & 1 \\
\end{bmatrix}\\
\alpha(a\circ w)
&=
\begin{bmatrix}
1 & 0 & 0 \\
1 & 1 & 0 \\
0 & 0 & 1 
\end{bmatrix}
\cdot
\begin{bmatrix}
x^{-x_w}t^{-t_w} & 0 & 0 \\
\num(w) & 1 & 0 \\
\den(w) & 0 & 1 \\
\end{bmatrix}
=
\begin{bmatrix}
x^{-x_w}t^{-t_w} & 0 & 0 \\
\num(w)+x^{-x_w}t^{-t_w} & 1 & 0 \\
\den(w) & 0 & 1 \\
\end{bmatrix}\\
\alpha(x^{-1}\circ w)
&=
\begin{bmatrix}
x & 0 & 0 \\
0 & 1 & 0 \\
x & 0 & 1 
\end{bmatrix}
\cdot
\begin{bmatrix}
x^{-x_w}t^{-t_w} & 0 & 0 \\
\num(w) & 1 & 0 \\
\den(w) & 0 & 1 \\
\end{bmatrix}
=
\begin{bmatrix}
x^{-(x_w -1)}t^{-t_w} & 0 & 0 \\
\num(w) & 1 & 0 \\
\den(w) + x^{-(x_w -1)}t^{-t_w} & 0 & 1 \\
\end{bmatrix}\\
\alpha(t^{-1}\circ w)
&=
\begin{bmatrix}
t & 0 & 0 \\
0 & 1 & 0 \\
0 & 0 & 1 \\
\end{bmatrix}
\cdot
\begin{bmatrix}
x^{-x_w}t^{-t_w} & 0 & 0 \\
\num(w) & 1 & 0 \\
\den(w) & 0 & 1 \\
\end{bmatrix}
=
\begin{bmatrix}
x^{-x_w}t^{-(t_w -1)} & 0 & 0 \\
\num(w) & 1 & 0 \\
\den(w) & 0 & 1 \\
\end{bmatrix}
\end{align*}}
In every case, the expressions in the first column coincide with the formulae from
Lemmas \ref{le:p-num-den1}, \ref{le:p-num-den3}, \ref{le:p-num-den5}, \ref{le:p-num-den2} and \ref{le:p-num-den4}.
Hence, the statement holds.
\end{proof}

\subsection{Tracing equation $w(a,t,x)$ on the $xt$-grid}
\label{se:tracing}

Define the $xt$-grid as the Cayley graph of 
$\MZ^2 = \Set{x^it^j}{i,j\in\MZ}$
in which every edge is labeled with $x^\pm,t^\pm$. 
Viewed this way, the $xt$-grid can be seen as an automaton with infinitely many states, 
with the initial state $x^0 t^0$.
The following statement is immediate.

\begin{lemma}\label{le:S-blf}
There is a one-to-one correspondence between 
$\delta$-parametric polynomials $\blf$ and 
finite subsets $S_{\blf} \subset \MZ^2$, given by
\[
\blf = \sum_{x^i t^j \in S_f} z^{i\delta + j}.
\]
\end{lemma}

For a word $w = w(a,t,x)$, define the finite subsets $N_w, D_w \subset \MZ^2$ by
\[
N_w = S_{\num(w)}, \qquad D_w = S_{\den(w)}.
\]
In this section, we describe a geometric procedure, called \emph{word tracing}, 
which constructs $N_w$ and $D_w$ by following the path on the $xt$-grid
determined by $w$, 
denoted $\Path(w)$, and adding points to $N_w$ and $D_w$ along the way.

Recall that the matrix formulae from the proof of 
Proposition~\ref{pr:num-den-matrices} compute the triple
\[
((-x_w, -t_w), \num_\delta(w), \den_\delta(w))
\]
by processing the word $w$ letter by letter from right to left.  
This computation can be interpreted as a geometric walk along a path $\Path(w)$, 
which is defined inductively by the following formulae 
(here, $\SD$ denotes the symmetric difference of sets).
\begin{itemize}
\item 
For the empty word $w = \varepsilon$, we have:  
\begin{itemize}
\item 
$\Path(\varepsilon)$ is the empty path, starting and ending at $x^0 t^0$,
\item 
the coordinates of the current point are $x_\varepsilon = t_\varepsilon = 0$,
\item
$N_\varepsilon = \varnothing$ and $D_\varepsilon = \varnothing$.
\end{itemize}
\item 
For $w'=x\circ w$
\begin{itemize}
\item 
$\Path(w')=\Path(w)\circ e$, where
$e=x^{-x_w}t^{-t_w}\to x^{-(x_w+1)}t^{-t_w}$;
\item 
$x_{w'}=x_{w}+1$ and $t_{w'}=t_{w}$;
\item 
$D_{w'} = D_w \SD \{x^{-x_w}t^{-t_w}\}$ and $N_{w'}=N_w$.
\end{itemize}
\item 
For $w'=x^{-1}\circ w$
\begin{itemize}
\item 
$\Path(w')=\Path(w)\circ e$,
where $e=x^{-x_w}t^{-t_w} \to x^{-(x_w-1)}t^{-t_w}$;
\item 
$x_{w'}=x_{w}-1$ and $t_{w'}=t_{w}$;
\item 
$D_{w'} = D_w \SD \{x^{-(x_w-1)}t^{-t_w}\}$ and $N_{w'}=N_w$.
\end{itemize}
\item 
For $w'=t\circ w$
\begin{itemize}
\item 
$\Path(w')=\Path(w)\circ e$,
where $e=x^{-x_w}t^{-t_w} \to x^{-x_w}t^{-(t_w+1)}$;
\item 
$x_{w'}=x_{w}$ and $t_{w'}=t_{w}+1$,
\item 
$D_{w'} = D_w$ and $N_{w'}=N_w$.
\end{itemize}
\item 
For $w'=t^{-1}\circ w$
\begin{itemize}
\item 
$\Path(w')=\Path(w)\circ e$,
where $e=x^{-x_w}t^{-t_w} \to x^{-x_w}t^{-(t_w-1)}$;
\item 
$x_{w'}=x_{w}$ and $t_{w'}=t_{w}-1$,
\item 
$D_{w'} = D_w$ and $N_{w'}=N_w$.
\end{itemize}
\item 
For $w'=a^{\pm}\circ w$
\begin{itemize}
\item 
$\Path(w')=\Path(w)$;
\item 
$x_{w'}=x_{w}$ and $t_{w'}=t_{w}$,
\item 
$D_{w'} = D_w$ and $N_{w'}=N_w \SD \{x^{-x_w}t^{-t_w}\}$.
\end{itemize}
\end{itemize}
To summarize, tracing the word $w$ letter by letter from right to left on the $xt$-grid proceeds as follows:  
\begin{itemize}
\item 
Reading $x$/$x^{-1}$: take a step left/right and add the current/right point to $D_w$ (mod 2), respectively.
\item
Reading $t$/$t^{-1}$: take a step down/up, respectively.
\item
Reading $a$/$a^{-1}$: remain at the current point and add it to $N_w$ (mod 2).
\end{itemize}

\begin{figure}[h]
\centering
\includegraphics[width=0.3\linewidth]{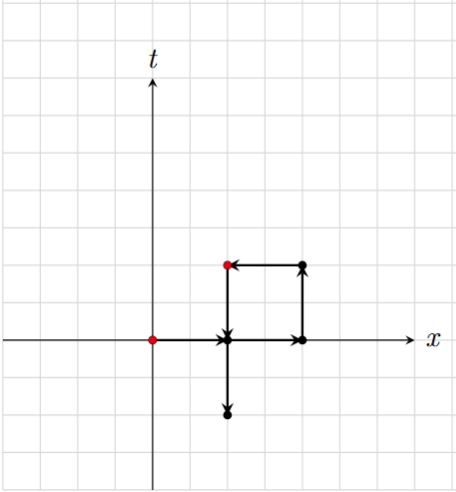}
\hspace{0.9cm}
\includegraphics[width=0.3\linewidth]{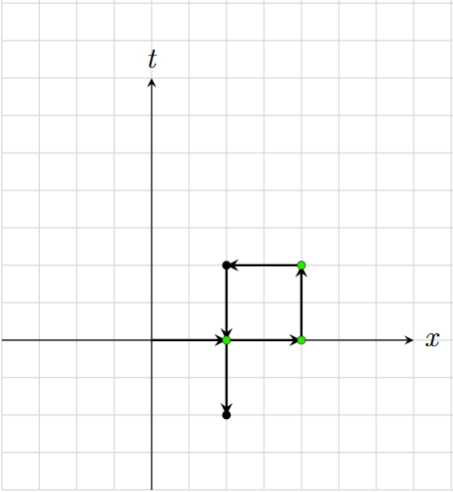}

\caption{
Tracing the word $w = t^2 a x t^{-1} x^{-2} a$ on the $xt$-grid produces the sets 
$N_w = \{ x^0 t^0, x^1 t^1 \}$ and $D_w = \{ x^1 t^0, x^2 t^0, x^2 t^1 \}$, 
corresponding to the polynomials $\num(w)$ and $\den(w)$, respectively.  
The left and right diagrams show $N_w$ (red points) and $D_w$ (green points), respectively.
}
\label{fig:tracing_example}
\end{figure}

\begin{lemma}\label{le:square}
For any word $w = w(a,t,x)$, the sets 
$N_w$ and $D_w$ are contained in the square 
$[-|w|, |w|] \times [-|w|, |w|]$.
\end{lemma}

\begin{proof}
Tracing a word $w$ of length $|w|$ on the $xt$-grid cannot move outside the square $[-|w|, |w|]^2$, 
so all points in $N_w$ and $D_w$ lie within this square.
\end{proof}

One also easily verifies that the sets $N_w$ and $D_w$
are contained in a $|w| \times |w|$ square of the integer grid.
Moreover, when $x_w = t_w = 0$, they are contained in a
$|w|/2 \times |w|/2$ square.


\subsection{Complexity of solving the equation $w(a,t,x)=1$ when $\sum \varepsilon_i \ne 0$}\label{subsection: complexity}

Suppose that the equation $w(a,t,x)=1$ satisfies $\sigma_x(w)\ne 0$.
Then the value of $\delta$ associated with the word $w(a,t,x)$ 
is defined by \eqref{eq:delta}.

\begin{lemma}\label{lemma:|w|}
$\sigma_x(w)\ne 0\ \ \Rightarrow\ \ |\delta| \le |w|$.
\end{lemma}

\begin{proof}
Because $|x_w|,|t_w|\le |w|$.
\end{proof}

Geometrically, $\delta = -\frac{t_w}{x_w}$ corresponds to 
the negative slope of the line 
passing through the origin $(0,0)$ and the terminus $(-x_w, -t_w)$ of $\Path(w)$.  
Substituting this value of $\delta$ 
into $\num(w)$ and $\den(w)$ 
is equivalent to 
projecting the points in $N_w$ and $D_w$ onto the $t$-axis 
along 
the direction of the vector $\vec{v} = (1, -\delta)$.

\begin{figure}[H]
\centering
\includegraphics[width=0.4\linewidth]{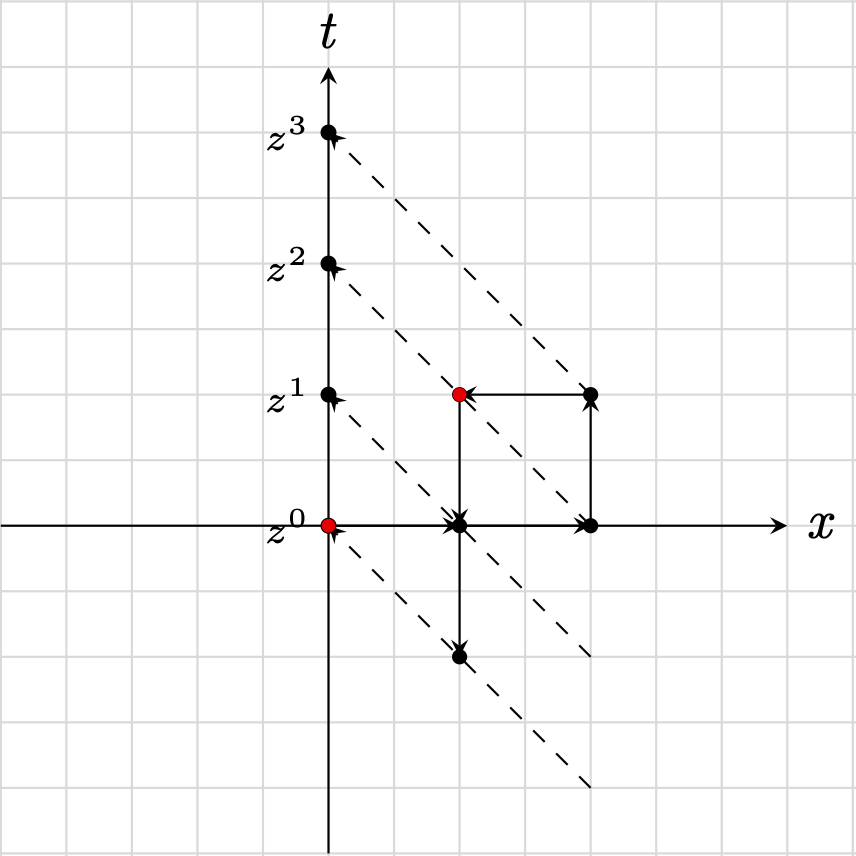}
\hspace{1cm}
\includegraphics[width=0.4\linewidth]{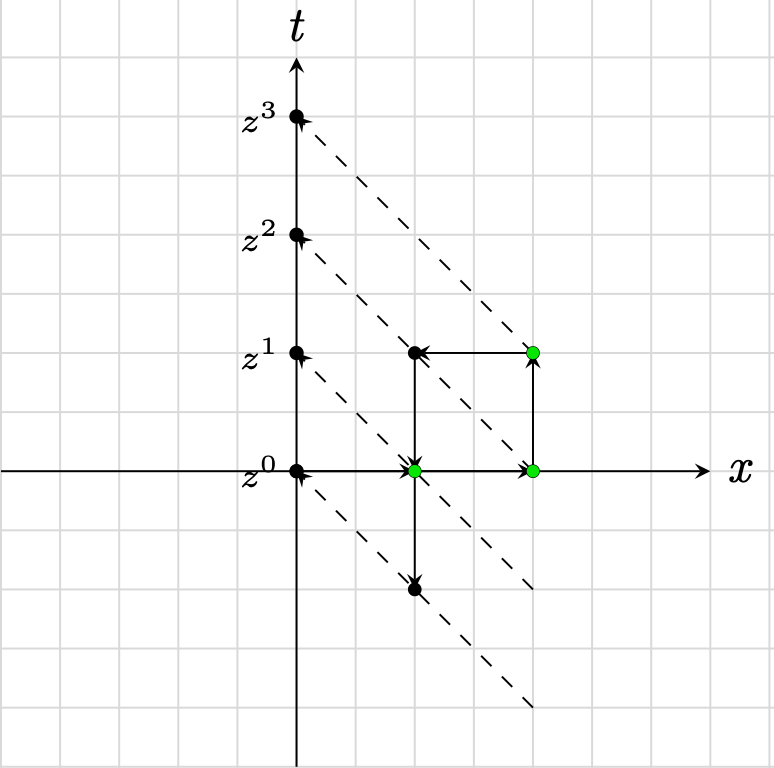}
\caption{Substituting $\delta \leftarrow 1$ to the example in Figure \ref{fig:tracing_example}, we can visualize the numerator (see left diagram) and denominator (see right diagram) polynomials on $t$-axis.
The sets $N_w$ and $D_w$ are $N_w = \{ z^0, z^2 \}$ and $D_w =\{ z^1, z^2, z^3 \}$. }
\label{fig:proj_onto_t}
\end{figure}

\begin{lemma}\label{le:trivial-den-num}
For any word $w = w(a,t,x)$,
$$
\num(w) \neq 0 \mbox{ and } |\delta| > |w|
\ \Rightarrow\ \num_\delta(w) \neq 0.
$$
An analogous statement holds for $\den(w)$.
\end{lemma}

\begin{proof}
Since $\num(w) \neq 0$, the set $N_w$ is nonempty.
If $|N_w| = 1$, then $\num_\delta(w) \neq 0$ for every $\delta \in \MZ$.
Assume therefore that $|N_w| > 1$.
Because $N_w$ is contained in a $|w| \times |w|$ square,
the absolute value of the slope of any non-vertical line determined by
two distinct points of $N_w$ is at most $|w|$.
Consequently, if $|\delta| > |w|$, then under the projection determined by
$\delta$ onto the $t$-axis, at least one point of $N_w$ survives,
and hence $\num_\delta(w) \neq 0$.
\end{proof}

\begin{corollary}\label{co:subst-deg-ord}
For any word $w=w(a,t,x)$, if $|\delta|\le |w|$ and
$\num_{\delta}(w)\ne 0$,
then the following holds:
$$
|\ord(\num_{\delta}(w))|,
|\deg(\num_{\delta}(w))|
\le |w|^2+|w|.$$
An analogous statement holds for $\den(w)$.
\end{corollary}

\begin{proof}
Projecting the point $(x,t)\in N_w$ onto the $t$-axis along the direction 
of the vector $\Vec{v} = (1, -\delta)$ creates the point
$(0,x\delta + t)$, which contributes the monomial
$z^{x\delta + t}$ to $\num_\delta(w)$. 
Hence, the square $[-|w|,|w|]^2$ of Lemma \ref{le:square}
is projected onto the interval
$[-|w|^2-|w|,|w|^2+|w|]$.
Therefore, nontrivial $\num_{\delta}(w)$ satisfies
$$
|\ord(\num_{\delta}(w))|,
|\deg(\num_{\delta}(w))|
\le |w|^2+|w|.
$$
The same proof works for $\den(w)$.
\end{proof}

The bounds of Corollary \ref{co:subst-deg-ord} are asymptotically sharp, 
as illustrated by the following example.

\begin{example}
Let us consider the word $w = t^{1-n} x^{n-1} a t^{-1} x^{-n}a$ of the length $|w| = 3n+ 1$. 
Its trace and the polynomials $\num(w)$ and $\den(w)$
are visualized in Figure \ref{fig:quadratic_den_and_num}.
Its terminus is the point $(1,n)$,
$x_w=-1$, $t_w=-n$, and $\delta = -n$.
Replacing $\delta$ with $-n$ in $\num(w)$ and $\den(w)$ we get
$\num_{\delta}(w) = z^{-n^2 + 1} + z^0$ and
$\den_{\delta}(w) = (z^{-n} + \dots + z^{-n^2}) + (z^{-2n+1} + \dots + z^{-n^2 + 1})
= \sum_{i=1}^n z^{-in} + \sum_{i=2}^n z^{-in + 1}$.

\begin{figure}[h]
    \centering
    \includegraphics[width=0.45\linewidth]{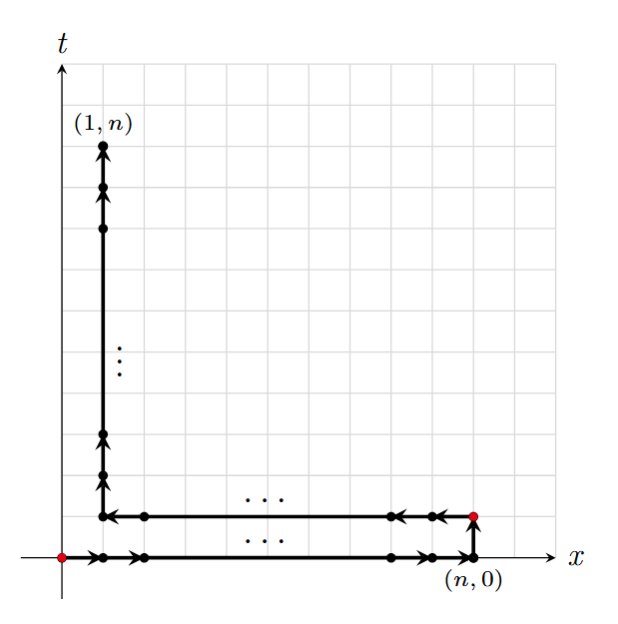}
    \hspace{0.5cm}
    \includegraphics[width=0.45\linewidth]{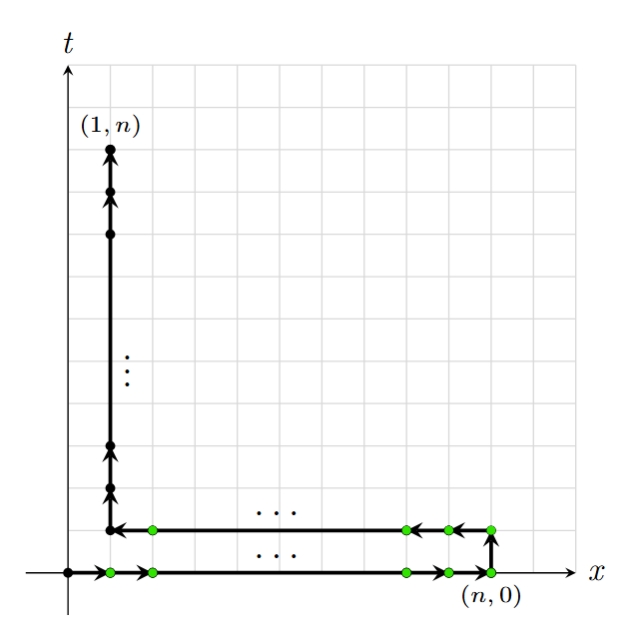}
\caption{Tracing the word $w = t^{1-n} x^{n-1} a t^{-1} x^{-n} a $ 
creates a path from $(0,0)$ to $(1,n)$ and produces the sets
$N_w =\{ x^n t^1, x^0 t^0 \}$ and 
$D_w =\{ x^1 t^0, x^2 t^0 \dots, x^n t^0, x^2 t^1, x^3 t^1, \dots , x^n t^1 \}$. 
}
    \label{fig:quadratic_den_and_num}
\end{figure}
\end{example}

\begin{proposition}\label{pr:compute-num-den}
If $|\delta| \le |w|$, then the representations 
$\pres(\num_\delta(w))$ and $\pres(\den_\delta(w))$ 
can be computed in $O(|w|^2)$ time.
\end{proposition}

\begin{proof}
The representation of $\num_\delta(w)$ is constructed in three steps.
First, allocate memory for the function
$n:[-|w|^2 - |w|, \, |w|^2 + |w|] \to \{0,1\}$ to store
the coefficients of $\num_\delta(w)$.
Next, trace the word $w$ as described 
in Section~\ref{se:tracing}.  
Instead of adding elements $x^{-x} t^{-t}$ to $N_w$, 
add $z^{x\delta + t}$ to $\num_\delta(w)$.  
Note that the multiplication by $\delta$ can be avoided: as the word is traced, 
it suffices to increment or decrement the exponent by $\delta$.  
Finally, process the resulting coefficients from left to right to form $\pres(\num_\delta(w))$.  
A similar procedure can be used to construct $\pres(\den_\delta(w))$.
\end{proof}

\begin{proposition}\label{Prop:algorithm-equations}
There exists an algorithm running in 
$O(|w|^2 \log |w| \log \log |w|)$ time that decides whether a given
one-variable equation $w(a,t,x) = 1$ with $\sigma_x(w) \neq 0$
admits a solution.
\end{proposition}

\begin{proof}
First, compute $x_w = \sigma_x(w)$ and $t_w = \sigma_t(w)$, and set
$\delta = -t_w / x_w$.
Division of integers of absolute value 
bounded by $|w|$ can be performed
in $O(|w|)$ time.
By Lemma~\ref{lemma:|w|}, we have $|\delta| \le |w|$.

Next, by Proposition~\ref{pr:compute-num-den},
the polynomials $\num_\delta(w)$ and $\den_\delta(w)$ 
can be computed in $O(|w|^2)$ time.
As in Proposition~\ref{pr:xw-nonzero}, it remains to check whether
\[
\den_\delta(w) = \num_\delta(w) = 0
\quad \text{or} \quad
\den_\delta(w) \mid \num_\delta(w).
\]
By Proposition~\ref{pr:Laurent-complexity}(d), this step requires
\[
O\!\left((|w|^2 + |w|)\,
\log(|w|^2 + |w|)\,
\log \log (|w|^2 + |w|)\right)
\]
time, which simplifies to
$O(|w|^2 \log |w| \log \log |w|)$.
\end{proof}

\section{Divisibility problem for parametric polynomials}
\label{se:divisibility-decidable}

Let $\blf$ and $\blg$ be $\delta$-parametric polynomials, with $\blf$
nontrivial. In this section we prove that 
$\DIV_+(\blf,\blg)$ is decidable.

\subsection{Normalization of $\blf$ and $\blg$}
Introduce a new variable $y$ and replace each occurrence of $z^\delta$
in $\blf$ and $\blg$ by $y$. This yields the bivariate polynomials
$F(z,y), G(z,y) \in \MZ_2[z^\pm,y^\pm]$
\begin{align*}
\blf&=\sum_{(i,j)\in S_{\scalebox{0.7}{$\blf$}}} z^{i\delta+j}
&\leftrightarrow&&
\sum_{(i,j)\in S_{\scalebox{0.7}{$\blf$}}} y^iz^j&=F(z,y)\\
\blg&=\sum_{(i,j)\in S_{\scalebox{0.7}{$\blg$}}} z^{i\delta+j}
&\leftrightarrow&&
\sum_{(i,j)\in S_{\scalebox{0.7}{$\blg$}}} y^iz^j&=G(z,y).
\end{align*}
Similarly to the definition of $S_\blf$ introduced in
Lemma~\ref{le:S-blf}, define the finite set
\[
S_F
=
\Set{x^i t^j \in \MZ^2}{y^i z^j \text{ occurs in } F},
\]
viewed as a subset of the $xt$-grid.
For a nontrivial $\delta$-parametric polynomial $\blf$ 
and the corresponding polynomial $F(z,y)$ define
\begin{align*}
\dmin_z(F) &= \min\Set{j}{y^i z^j\text{ is involved in }F},
&\dmin_z(\blf) &= \min\Set{j}{z^{\,i\delta+j}\text{ is involved in }\blf },\\
\dmax_z(F) &= \max\Set{j}{y^i z^j\text{ is involved in }F},
&\dmax_z(\blf) &= \max\Set{j}{z^{\,i\delta+j}\text{ is involved in }\blf },\\
\dmin_y(F) &= \min\Set{i}{y^i z^j\text{ is involved in }F},
&\dmin_y(\blf) &= \min\Set{i}{z^{\,i\delta+j}\text{ is involved in }\blf },\\
\dmax_y(F) &= \max\Set{i}{y^i z^j\text{ is involved in }F},
&\dmax_y(\blf) &= \max\Set{i}{z^{\,i\delta+j}\text{ is involved in }\blf }.
\end{align*}
These characteristics are defined as $\varnothing$ when $\blf=0$ (or $F=0$).

Observe that  multiplying $\blf$ and $\blg$
by arbitrary powers of $z$ and $z^\delta$ does not affect the value 
of $\DIV_+(\blf,\blg)$:
\[
\DIV_+(\blf,\blg)
=
\DIV_+(\blf\, z^{i\delta+i'},\blg\, z^{j\delta+j'}),\quad
\text{for all }i,j,i',j'\in\MZ.
\]
Thus, without loss of generality, for nontrivial $\blf,\blg$
we may assume that
\[
\dmin_z(\blf) = \dmin_y(\blf) =
\dmin_z(\blg) = \dmin_y(\blg) = 0.
\]
which implies that
\[
\blf=\sum_{i=0}^s f_i(z) z^{i\delta},\qquad
\blg=\sum_{i=0}^t g_i(z) z^{i\delta},
\]
where $f_0,g_0 \ne 0$ and $f_i,g_i \in \MZ_2[z]$ are ordinary
polynomials (rather than Laurent polynomials).
We refer to this step as the \emph{normalization of monomials}
in $\blf$ and $\blg$ (similarly, in $F$ and $G$).

\begin{lemma}\label{le:S_F-bounds}
$S_{\scalebox{0.7}{$\blf$}} = S_F \subseteq 
[0,s]\times[0,\dmax_z(\blf)]$, where $s = \dmax_y(\blf)$.
\end{lemma}

\begin{proof}
The equality $S_{\scalebox{0.7}{$\blf$}} = S_F$ 
follows directly from the construction of the polynomial $F(z,y)$.
The inclusion
$S_F \subseteq [0,\dmax_y(\blf)]\times[0,\dmax_z(\blf)]$
follows from the normalization of monomials in $\blf$.
\end{proof}

\subsection{Dividing $G$ by $F$}

Treat $F$ and $G$ as polynomials in the variable $y$ with coefficients
in the fraction field $\MZ_2(z)$ of $\MZ_2[z]$, and perform Euclidean division of $G$ by $F$.
This yields
\begin{equation}\label{eq:G-divide-F}
G(z,y) = F(z,y)\, q(z,y) + r(z,y),
\end{equation}
where $q(z,y), r(z,y) \in \MZ_2(z)[y]$ and either $r=0$ identically or 
$\deg_y(r) < \deg_y(F)$.

\subsubsection{Case when $t=s$}
The quotient and remainder $q(z,y)$ and $r(z,y)$ in~\eqref{eq:G-divide-F}
can be computed explicitly by successive elimination of the leading
monomials of $G$ using suitable multiples of $F$.
To illustrate this process, consider a single elimination step, namely,
the case when $F$ and $G$ have the same degree in $y$,
that is when $t=s$:
\begin{align*}
F(z,y) &= f_s(z)\, y^s + \dots + f_1(z)\, y + f_0(z),\\
G(z,y) &= g_s(z)\, y^s + \dots + g_1(z)\, y + g_0(z).
\end{align*}
In this case, a single division step yields
\begin{equation}\label{eq:rq-one-step}
q(z,y) = \frac{g_s(z)}{f_s(z)},
\quad
r(z,y) = \sum_{i=0}^{s-1}
\left(g_i(z) - \frac{g_s(z)}{f_s(z)}\, f_i(z)\right) y^i.
\end{equation}

The computation of the remainder $r(z,y)$ can be expressed in vector form
as an affine linear transformation of the coefficient vector of $G$.
More precisely, we have
\[
\begin{pmatrix}
g_s\\
g_{s-1}\\
\vdots\\
g_2\\
g_1
\end{pmatrix}
\;\stackrel{g_0}{\longmapsto}\;
\begin{pmatrix}
g_{s-1} - g_s \tfrac{f_{s-1}}{f_s}\\
g_{s-2} - g_s \tfrac{f_{s-2}}{f_s}\\
\vdots\\
g_{1} - g_s \tfrac{f_{1}}{f_s}\\
g_{0} - g_s \tfrac{f_{0}}{f_s}
\end{pmatrix}
=
\frac{1}{f_s}\, M_f
\begin{pmatrix}
g_s\\
g_{s-1}\\
\vdots\\
g_2\\
g_1
\end{pmatrix}
+
g_0\,\ove_s,
\]
where
\[
M_f =
\begin{pmatrix}
-f_{s-1} & f_s & 0 & \cdots & 0\\
-f_{s-2} & 0 & f_s & \cdots & 0\\
\vdots & & & \ddots & \\
-f_{1} & 0 & 0 & \cdots & f_s\\
-f_{0} & 0 & 0 & \cdots & 0
\end{pmatrix},
\qquad
\ove_s =
\begin{pmatrix}
0\\
0\\
\vdots\\
0\\
1
\end{pmatrix}.
\]

\subsubsection{Case when $t>s$}
In this case the remainder $r(z,y)$ is obtained by
iterating the above transformation. 
Starting from the initial coefficient vector
\[
\ovu = \ovu_0 =
\begin{pmatrix}
g_t\\
g_{t-1}\\
\vdots\\
g_{t-s+1}
\end{pmatrix},
\]
we obtain a sequence of vectors
\[
\ovu_0
\;\stackrel{g_{t-s}}{\longmapsto}\;
\ovu_{1}
\;\stackrel{g_{t-s-1}}{\longmapsto}\;
\ovu_{2}
\;\stackrel{g_{t-s-2}}{\longmapsto}\;
\;\cdots\;
\stackrel{g_{0}}{\longmapsto}\;
\ovu_{t-s+1},
\]
where
\[
\ovu_{k+1}
=
\frac{1}{f_s} M_f \ovu_k + g_{t-s-k}\,\ove_s,
\qquad k=0,1,\dots,t-s.
\]
The final vector $\ovu_{t-s+1}$ encodes the coefficients of the
remainder $r(z,y)$ and admits the explicit expression
\[
\ovu_{t-s+1}
=
\left(\frac{1}{f_s} M_f\right)^{t-s+1}\ovu
+
\sum_{k=0}^{t-s}
\left(\frac{1}{f_s} M_f\right)^{k}\, g_k\,\ove_s.
\]

\begin{lemma}\label{le:u_ij}
For every $i=0,\dots,t-s+1$, the following statements hold:
\begin{enumerate}
\item[(a)]
The vector $\ovu_i$ can be written in the form
$\ovu_i = \frac{1}{f_s^i}(u_{i1},\dots,u_{is})^T$,
where $u_{ij}\in\MZ_2[z]$ are ordinary polynomials.
\item[(b)]
$\deg(u_{ij})
\le
\dmax_z(\blg)+i\,\dmax_z(\blf)$
for all $j=1,\dots,s$.
\item[(c)]
The vectors $\ovu_1,\dots,\ovu_{t-s+1}$ can be computed in 
\begin{equation}\label{eq:u_ij-complexity}
\wO \bigl( \dmax_y(\blg)\dmax_y(\blf)
(\dmax_z(\blg) + \dmax_y(\blg)\dmax_z(\blf)) \bigr)
\end{equation}
time.
\end{enumerate}
\end{lemma}

\begin{proof}
We prove (a) by induction on $i$.
For $i=0$, the claim holds with $u_{0j}=g_{t-j+1}$.
Assume that (a) holds for some $i\in\MN$. Then
\[
\ovu_{i+1}
=
\frac{1}{f_s}M_f\ovu_i+g_{t-s-i}\,\ove_s
=
\frac{1}{f_s^{i+1}}
M_f(u_{i1},\dots,u_{is})^T
+
g_{t-s-i}\,\ove_s.
\]
Since the entries of $M_f$ lie in $\MZ_2[z]$, it follows that
\[
\ovu_{i+1}
=
\frac{1}{f_s^{i+1}}
(u_{i+1,1},\dots,u_{i+1,s})^T
\]
for some polynomials $u_{i+1,j}\in\MZ_2[z]$, proving~(a).

The transition from $\ovu_i$ to $\ovu_{i+1}$ can be written explicitly as
\[
\ovu_i=
\begin{pmatrix}
u_{i,s}\\
u_{i,s-1}\\
\vdots\\
u_{i,2}\\
u_{i,1}
\end{pmatrix}
\stackrel{g_{t-s-i}}{\longmapsto}
f_s
\begin{pmatrix}
u_{i,s-1}\\
u_{i,s-2}\\
\vdots\\
u_{i,1}\\
g_{t-s-i}
\end{pmatrix}
-
u_{i,s}
\begin{pmatrix}
f_{s-1}\\
f_{s-2}\\
\vdots\\
f_{1}\\
f_{0}
\end{pmatrix}
=
\ovu_{i+1}.
\]
Each transition increases the degree of every entry by at most
$\dmax_z(\blf)$. Therefore,
\[
\deg(u_{i+1,j})
\le
\dmax_z(\blg)+(i+1)\dmax_z(\blf),
\]
proving~(b).

Computation of all $u_{ij}$ involves
$2s(t-s+1)$ multiplications and $t-s+1$ additions
of polynomials whose degrees satisfy
\begin{align*}
\deg(u_{ij})
&\le
\dmax_z(\blg)+(t-s+1)\dmax_z(\blf)\\
&\le
\dmax_z(\blg)+\dmax_y(\blg)\dmax_z(\blf).
\end{align*}
Hence, the total running time is
\[
\wO\!\Bigl(
\dmax_y(\blg)\dmax_y(\blf)
\bigl(
\dmax_z(\blg)
+
\dmax_y(\blg)\dmax_z(\blf)
\bigr)
\Bigr).
\eqno\qed
\]\let\qed\relax
\end{proof}

\begin{lemma}\label{le:Q-R-rational}
There exist polynomials $Q(z,y),R(z,y)\in\MZ_2[z,y]$ satisfying
\[
q(z,y) = \frac{Q(z,y)}{f_s^{t-s+1}},
\quad
r(z,y) = \frac{R(z,y)}{f_s^{t-s+1}}.
\]
Moreover, the polynomial $R(z,y)$ can be computed in time
\eqref{eq:u_ij-complexity}.
\end{lemma}

\begin{proof}
The coefficients of $r(z,y)$ are determined by the vector
$\ovu_{t-s+1}$.
By Lemma~\ref{le:u_ij}, the entries of $\ovu_{t-s+1}$ have the form
$\frac{u_{t-s+1,j}}{f_s^{\,t-s+1}}$,
where $u_{t-s+1,j}\in\MZ_2[z]$. By
Lemma~\ref{le:u_ij}(c), these polynomials can be computed in time
\eqref{eq:u_ij-complexity}.
Hence the coefficients of $R(z,y)$ can be computed within the same
time bound.

The coefficients of $q(z,y)$ are produced iteratively during the long division process.
At the $i$th step, by \eqref{eq:rq-one-step}, the newly generated coefficient is
\[
\frac{\frac{u_{i-1,i}}{f_s^{\,i-1}}}{f_s} = \frac{u_{i-1,i}}{f_s^{\,i}},
\]
and therefore has the required form.
Collecting these coefficients yields the desired representation of
$q(z,y)$.
\end{proof}

\begin{lemma}\label{le:S_R-bounds}
$S_R \subseteq [0,s-1]\times [0,\dmax_z(\blg)+(t-s+1)\dmax_z(\blf)]$.
\end{lemma}

\begin{proof}
Follows from Lemma~\ref{le:u_ij}(b).
\end{proof}

By definition of $\DIV_+$ and by construction of the polynomials
$F(z,y)$ and $G(z,y)$, we have
\begin{align*}
\DIV_+(\blf,\blg)=true
&\quad\Leftrightarrow\quad
\exists\,\delta>0
\text{ such that }
\blf_\delta\mid\blg_\delta\\
&\quad\Leftrightarrow\quad
\exists\,\delta>0
\text{ such that }
F(z,z^\delta)\mid G(z,z^\delta),
\end{align*}
because 
$\blf_\delta = F(z, z^\delta)$ and 
$\blg_\delta = G(z, z^\delta)$.
To determine whether $F(z, z^\delta) \mid G(z, z^\delta)$
we distinguish two cases: $r(z,y)=0$ and $r(z,y)\ne 0$.

\subsubsection{Case when $t<s$}

In the case $t<s$, division of $G$ by $F$ involves no elimination steps,
and \eqref{eq:G-divide-F} takes the form
\[
G(z,y)=F(z,y)\cdot 0+G(z,y).
\]
If $G=0$, then $\DIV_+(\blf,\blg)=true$, and there is nothing further
to prove. Otherwise, if $G\ne 0$, then the remainder of division is
$r=R=G\ne 0$,
which falls into the case considered in Section~\ref{se:r-ne-0}.
In this situation, we may set $t-s+1=0$
and all subsequent estimates, including the polynomial-time complexity
bounds, remain valid.

In what follows we assume that $t\ge s$.

\subsection{Case when $r(z,y)=0$}
\label{se:r-0}

Before discussing the case \(r(z,y)=0\), let us recall Gauss's lemma.
Let $R$ be a Euclidean domain.
For a nonzero polynomial
$f(y)=\alpha_p y^p+\dots+\alpha_1 y+\alpha_0 \in R[y]$,
the \emph{content} of $f$, denoted by $\cont(f)$, is defined as
$$
\cont(f)=\gcd(\alpha_p,\dots,\alpha_0)\in R,
$$
where the $\gcd$ is chosen up to multiplication by a unit of $R$.
A polynomial $f\in R[y]$ is called \emph{primitive} if
$\cont(f)=1$.
Every nonzero polynomial $f\in R[y]$ can be written in the form
$f=\cont(f)\cdot f_0$,
where $f_0\in R[y]$ is primitive.
By the classical Gauss lemma, the product of two primitive polynomials 
is primitive.

Consider the Euclidean domain $R=\MZ_2[z]$ whose only unit is $1$,
and the polynomials $F(z,y),G(z,y),Q(z,y)\in \MZ_2[z][y]$
with coefficients from $R$.

\begin{lemma}
For a normalized polynomial $F$, the content
$\cont(F)$ can be computed in time $O(\dmax_y(F)\dmax_z^2(F))$.
\end{lemma}

\begin{proof}
We have $S_F\subseteq [0,\deg_y(F)] \times [0,\deg_z(F)]$.
Computing $\cont(F)$ amounts to applying the Euclidean algorithm
to the coefficient polynomials of $F$.
Since each coefficient polynomial is represented by a $0/1$ column
determined by $S_F$, the total number of elementary steps
can be bounded by $\dmax_y(F)\dmax_z^2(F)$
yielding the stated complexity bound.
\end{proof}

Write 
\[
F=\cont(F)\cdot F_0,\qquad
G=\cont(G)\cdot G_0,\qquad
Q=\cont(Q)\cdot Q_0,
\]
for the primitive decompositions of $F,G,Q$.
Now, the following holds:
\begin{align*}
r(z,y)=0
&\ \stackrel{\scalebox{0.6}{\eqref{eq:G-divide-F}}}{\Leftrightarrow}\ 
G(z,y) = F(z,y)\cdot \frac{Q(z,y)}{f_s^{t-s+1}(z)} \\
&\ \Leftrightarrow\ 
G_0(z,y) = F_0(z,y)Q_0(z,y)\cdot \frac{\cont(F)\cont(Q)}{f_s^{t-s+1}(z)\cont(G)}\\
&\ \Leftrightarrow\
G_0(z,y)=u\,F_0(z,y)Q_0(z,y)
\text{ for a unit }u\in R\\
&\ \Leftrightarrow\ 
G_0(z,y) = F_0(z,y)Q_0(z,y) \\
&\ \Leftrightarrow\ 
G(z,y) \cont(F) = \cont(G) F(z,y)Q_0(z,y)\\
&\ \Leftrightarrow\ 
G(z,y) = F(z,y)Q_0(z,y)\frac{\cont(G)}{\cont(F)}.
\end{align*}
It follows that, when $r(z,y)=0$,
\begin{align*}
\DIV_+(\blf,\blg)=true
&\quad\Leftrightarrow\quad
\exists\,\delta>0
\text{ such that }
F(z,z^\delta)\mid G(z,z^\delta)\\
&\quad\Leftrightarrow\quad
\exists\,\delta>0
\text{ such that }
Q_0(z,z^\delta)\cont(G) \equiv 0 \bmod \cont(F).
\end{align*}

\begin{proposition}\label{pr:delta-bound1}
Let $\blf$ and $\blg$ be $\delta$-parametric polynomials, with $\blf$
nontrivial, and assume that $r(z,y)=0$.
Then $\DIV_+(\blf,\blg) =$ true
if and only if
\[
Q_0(z,z^\delta)\cont(G) \equiv 0 \bmod \cont(F)
\]
for some integer $\delta$ satisfying 
$1 \le \delta \le 2^{\deg(\cont(F))}$.
\end{proposition}

\begin{proof}
Let $N=2^{\deg(\cont(F))}$ and consider the sequence of remainders
\[
z^\delta \bmod \cont(F),\qquad \delta=1,2,\dots.
\]
There are at most $N$ distinct remainders modulo $\cont(F)$.
Hence some remainder repeats among the first $N+1$ terms.
Once a remainder repeats, all subsequent remainders repeat periodically.
Therefore every remainder appearing in the sequence appears for some
$1\le \delta\le N+1$.

Since polynomial evaluation modulo $\cont(F)$ depends only on the residue
class of $z^\delta$ modulo $\cont(F)$, the same conclusion holds for
\[
Q_0(z,z^\delta)\cont(G) \bmod \cont(F).
\]
Thus, if
\[
Q_0(z,z^\delta)\cont(G)\equiv 0 \pmod{\cont(F)}
\]
holds for some $\delta>0$, then it holds for some
$1\le \delta\le 2^{\deg(\cont(F))}$.
The claim follows.
\end{proof}

\begin{lemma}\label{le:verify-delta1}
If $r(z,y)=0$, then there is an algorithm that, for a given integer
$1\le \delta\le 2^{\deg(\cont(F))}$ 
determines whether
\[
Q_0(z,z^\delta) \cont(G)\equiv 0 \pmod{\cont(F)}
\]
in polynomial time in $\deg(\cont(F))$.
\end{lemma}

\begin{proof}
Using repeated squaring, one can compute the sequence of remainders
\[
\{z^{2^i} \bmod \cont(F)\}_{i=0}^{\deg(\cont(F))}.
\]
Express $\delta$ in binary form,
$\delta=\sum_{i=0}^m \varepsilon_i 2^i$,
where $\varepsilon_i\in\{0,1\}$,
and use the identity
\[
z^\delta
=
\prod_{\varepsilon_i=1} z^{2^i}
\]
to compute $z^\delta \bmod \cont(F)$.
Substituting this value into $Q_0(z,y) \cont(G)$ yields
\[
Q_0(z,z^\delta) \cont(G)\bmod \cont(F).
\]
All computations involve polynomial arithmetic modulo
$\cont(F)$ with bounded degrees, and therefore
require polynomial time.
\end{proof}

\begin{proposition}\label{pr:r=0-complexity}
If $r(z,y)=0$, then $\DIV_+(\blf,\blg)$ can be decided in time
$2^{O(\deg(\cont(F)))}$.
\end{proposition}

\begin{proof}
The algorithm checks, for each integer
$\delta=1,\dots,2^{\deg(\cont(F))}$ 
whether 
$Q_0(z,z^\delta)\cont(G) \equiv 0 \pmod{\cont(F)}$.
If the congruence holds for some $\delta$, then
$\DIV_+(\blf,\blg)=true$; otherwise,
$\DIV_+(\blf,\blg)=false$.

By Lemma~\ref{le:verify-delta1}, each individual value of $\delta$
can be processed in time polynomial in
$\deg(\cont(F))$.
Since the algorithm performs at most
$2^{\deg(\cont(F))}$ such checks, the total running time is
$2^{O(\deg(\cont(F)))}$.
\end{proof}

\subsection{Case when $r(z,y)\ne 0$}
\label{se:r-ne-0}

If $r(z,y)\ne 0$, then, by Lemma~\ref{le:Q-R-rational},
equation~\eqref{eq:G-divide-F} can be rewritten as
\[
R(z,y) = f_s^{\,t-s+1}G(z,y)-F(z,y)\,Q(z,y).
\]
Consequently,
\begin{align*}
\DIV_+(\blf,\blg)=true
&\quad\Leftrightarrow\quad
\exists\,\delta>0
\text{ such that }
F(z,z^\delta)\mid G(z,z^\delta)\\
&\quad\Rightarrow\quad
\exists\,\delta>0
\text{ such that }
F(z,z^\delta)\mid R(z,z^\delta).
\end{align*}

\begin{proposition}\label{pr:delta-bound2}
Suppose that $r(z,y)\ne 0$ and let $\delta>0$.
If $F(z, z^\delta) \mid G(z, z^\delta)$, then
$\delta \le \dmax_z(\blg)+(t-s+1)\dmax_z(\blf)$.
\end{proposition}

\begin{proof}
By Lemmas~\ref{le:S_F-bounds} and~\ref{le:S_R-bounds}, we have
\begin{align*}
S_F
&\subseteq
[0,s]\times[0,\dmax_z(\blf)],\\
S_R
&\subseteq
[0,s-1]\times
[0,\dmax_z(\blg)+(t-s+1)\dmax_z(\blf)].
\end{align*}
If $\delta > \dmax_z(\blf)$, then the monomials arising from the terms
$z^i y^s$ in $F(z,y)$ do not cancel after the substitution
$y=z^\delta$.
Consequently, if $\delta > \dmax_z(\blf)$, then
\[
s\delta
\le
\deg\!\bigl(F(z,z^\delta)\bigr)
\le
s\delta + \dmax_z(\blf).
\]
Similarly, if
$\delta > \dmax_z(\blg)+(t-s+1)\dmax_z(\blf)$,
then the monomials arising from the terms $z^i y^{s-1}$ in $R(z,y)$
do not cancel with monomials from lower $y$-degree terms after the
substitution $y=z^\delta$.
Consequently, if $\delta > \dmax_z(\blg)+(t-s+1)\dmax_z(\blf)$, then 
\[
(s-1)\delta
\le
\deg\!\bigl(R(z,z^\delta)\bigr)
\le
(s-1)\delta
+
\dmax_z(\blg)+(t-s+1)\dmax_z(\blf).
\]
Therefore, if $\delta > \dmax_z(\blg)+(t-s+1)\dmax_z(\blf)$, then the following holds:
\begin{align*}
F(z,z^\delta) \mid G(z,z^\delta)
&\ \Rightarrow\ 
\deg\!\bigl(F(z,z^\delta)\bigr)
\le
\deg\!\bigl(G(z,z^\delta)\bigr) \\
&\ \Rightarrow\ 
s\delta
\le
(s-1)\delta
+
\dmax_z(\blg)
+
(t-s+1)\dmax_z(\blf) \\
&\ \Rightarrow\ 
\delta \le \dmax_z(\blg)+(t-s+1)\dmax_z(\blf).
\end{align*}
The obtained contradiction shows that 
$\delta > \dmax_z(\blg)+(t-s+1)\dmax_z(\blf)$ is impossible.
\end{proof}

\begin{proposition}\label{pr:r-not0-complexity}
If $r(z,y)\ne 0$ for normalized $\blf$ and $\blg$, 
then $\DIV_+(\blf,\blg)$ can be decided in time
\begin{equation}\label{eq:r-not0-complexity}
\wO\!\rb{
(\dmax_y(\blf)+1)
\rb{
\dmax_z(\blg)
+
(t-s+1)\dmax_z(\blf)
}^2}.
\end{equation}
\end{proposition}

\begin{proof}
By Proposition~\ref{pr:delta-bound2}, it suffices to consider integers
\[
\delta=1,\dots,\dmax_z(\blg)+(t-s+1)\dmax_z(\blf).
\]
For each such $\delta$, construct the polynomials
$\blf_\delta$ and $\blg_\delta$.
Their degrees satisfy
\begin{align*}
\deg(\blf_\delta),\deg(\blg_\delta) 
&\le 
s(\dmax_z(\blg)+(t-s+1)\dmax_z(\blf)) + \max(\dmax_z(\blf),\dmax_z(\blg))\\
&\le 
(s+1)(\dmax_z(\blg)+(t-s+1)\dmax_z(\blf)).
\end{align*}
For each pair of polynomials $\blf_\delta,\blg_\delta$, testing whether
$\blf_\delta\mid\blg_\delta$
requires time
$$
\wO((s+1)(\dmax_z(\blg)+(t-s+1)\dmax_z(\blf)))
$$
Multiplying by the number of values of $\delta$ yields the claimed
running time.
\end{proof}

\subsection{Divisibility problem for parametric polynomials is decidable}

First, we prove that $\DIV_+$ is decidable and, moreover, belongs to $\NP$.

\begin{theorem}\label{th:DIV+NP}
$\DIV_+\in\NP$.
\end{theorem}

\begin{proof}
Let $\blf$ and $\blg$ be $\delta$-parametric polynomials.
In polynomial time, we compute the associated polynomials
$F(z,y)$ and $G(z,y)$, together with the quotient and remainder
$q(z,y)$ and $r(z,y)$ obtained by division of $G$ by $F$
in $\MZ_2(z)[y]$.
By construction, $(\blf,\blg)$ is a positive instance of $\DIV_+$
if and only if $F(z,z^\delta) \mid G(z,z^\delta)$ for some $\delta>0$.

\medskip\noindent
\textsc{Case 1.}
If $r(z,y)=0$, then, 
by Lemma~\ref{pr:delta-bound1} we may assume that
$\delta \le 2^{\deg(\cont(F))}$ and
by Lemma~\ref{le:verify-delta1},
the condition 
$F(z,z^\delta) \mid G(z,z^\delta)$ can be verified in polynomial time.

\medskip\noindent
\textsc{Case 2.}
If $r(z,y)\ne 0$, then, 
by Lemma~\ref{pr:delta-bound2} we may assume that
$\delta \le \dmax_z(\blg)+(t-s+1)\dmax_z(\blf)$, in which case
we directly check if $\blf_\delta \mid \blg_\delta$.

\medskip\noindent
Thus, in both cases, there exists a value of $\delta$
whose binary encoding has polynomial size and whose validity
can be verified in polynomial time.
\end{proof}

By Proposition~\ref{pr:div-div2}, divisibility for $(\blf,\blg)$
reduces to positive divisibility:
\[
\DIV(\blf,\blg)
=
\DIV_+(\blf,\blg)
\;\vee\;
(\blf_0\mid\blg_0)
\;\vee\;
\DIV_+(\rev(\blf),\rev(\blg)).
\]
We therefore next compare the relevant characteristics of
$(\blf,\blg)$ and $(\rev(\blf),\rev(\blg))$.

In Section~\ref{se:DIV-problem}, for every $\delta$-parametric polynomial
$\blf$ we defined the $\delta$-parametric polynomial $\rev(\blf)$.
Similarly, we introduce an analogous map for bivariate Laurent polynomials:
\[
\rev:\MZ_2(z)[y^\pm]\to \MZ_2(z)[y^\pm],
\]
defined by fixing the coefficients from $\MZ_2(z)$ and sending
$y\mapsto y^{-1}$.
The map $\rev$ is a ring automorphism.

As above, let $\blf,\blg$ be normalized
$\delta$-parametric polynomials, and let $F,G \in \MZ_2(z)[y^\pm]$
be the corresponding polynomials. 
Consider also the (normalized) $\delta$-parametric polynomials
$\rev(\blf),\rev(\blg)$ 
and their corresponding polynomials
$F',G' \in \MZ_2(z)[y^\pm]$.

\begin{lemma}\label{le:rev-characteristics}
The following holds:
\begin{itemize}
\item
$\cont(F)=\cont(F')$ and $\cont(G)=\cont(G')$.
\item
$\dmin_y(F)=\dmin_y(F')=\dmin_z(F)=\dmin_z(F')=0$.
\item
$\dmin_y(G)=\dmin_y(G')=\dmin_z(G)=\dmin_z(G')=0$.
\item
$\dmax_y(F)=\dmax_y(F')$ and $\dmax_z(F)=\dmax_z(F')$.
\item
$\dmax_y(G)=\dmax_y(G')$ and $\dmax_z(G)=\dmax_z(G')$.
\end{itemize}
\end{lemma}

\begin{proof}
The claim follows directly from the definition of $\rev$ and from
the normalization procedure.
\end{proof}

Divide $G$ by $F$ to obtain the remainder
$r$, and divide $G'$ by $F'$ to obtain the remainder $r'$.

\begin{lemma}\label{le:r=r'=0}
$r=0\ \Leftrightarrow\ r'=0$.
\end{lemma}

\begin{proof}
($\Rightarrow$)
Suppose that $r=0$, that is $F\mid G$. Then
$$
G(z,y) = F(z,y)\, q(z,y),
$$
for some $q(z,y) \in \MZ_2(z)[y]$. Then
$$
\rev(G(z,y)) = \rev(F(z,y))\, \rev(q(z,y)).
$$
By construction, $F'$ and $G'$ are obtained from
$\rev(F)$ and $\rev(G)$ by normalization of monomials.
Since normalization only multiplies polynomials by powers of $y$ and $z$,
it preserves divisibility.
Hence, $F'\mid G'$ and $r'=0$.

\medskip\noindent
($\Leftarrow$)
The converse is proved similarly.
\end{proof}

\begin{theorem}\label{th:DIV-NP}
The problem of determining whether a given pair $(\blf,\blg)$
is a positive instance of $\DIV$ is decidable and belongs to $\NP$.
\end{theorem}

\begin{proof}
By Theorem~\ref{th:DIV+NP}, the problem $\DIV_+(\blf,\blg)$
belongs to $\NP$.
By Lemma~\ref{le:rev-characteristics},
the pair $(\rev(\blf),\rev(\blg))$ has the same relevant
characteristics as $(\blf,\blg)$.
Hence,
$\DIV_+(\rev(\blf),\rev(\blg))$
also belongs to $\NP$.
Finally, by Proposition~\ref{pr:div-div2},
\[
\DIV(\blf,\blg)
=
\DIV_+(\blf,\blg)
\;\vee\;
(\blf_0\mid\blg_0)
\;\vee\;
\DIV_+(\rev(\blf),\rev(\blg)).
\]
Therefore, $\DIV(\blf,\blg)$ is decidable and belongs to $\NP$.
\end{proof}

\subsection{One-variable equations over the lamplighter group are decidable}
\label{se:equations-complexity}

We now apply the results obtained for the divisibility problem to
establish decidability and analyze the complexity of solving
one-variable equations over the lamplighter group,
distinguishing several cases.

\begin{lemma}\label{le:delta-bound2}
Let $w(a,t,x)=1$ be a one-variable equation, and let
$\blf$ and $\blg$ be the normalized forms of $\num(w)$ and $\den(w)$,
respectively. Then
\begin{itemize}
\item[(a)]
$0\le
\dmax_z(\blg),\dmax_z(\blf),
\dmax_y(\blg),\dmax_y(\blf)
\le |w|$.
\item[(b)]
$0 \le t-s+1 \le |w|$, provided that  $t\ge s$.
\item[(c)]
$\deg(\cont(F))\le |w|$,
where $F$ is the bivariate polynomial corresponding to $\blf$.
\end{itemize}
\end{lemma}

\begin{proof}
The first inequality follows directly from normalization.
The second inequality holds because
\[
t-s = \dmax_y(\blg)-\dmax_y(\blf) \le |w|.
\]
Since $\deg(f_i)\le |w|$ for every coefficient polynomial $f_i$,
we have
\[
\deg(\cont(F))
\le
\min_i\deg(f_i)
\le |w|,
\]
proving the last claim.
\end{proof}

\begin{theorem}\label{th:one-var-complexity}
The problem of determining whether a given one-variable equation
$w(a,t,x)=1$ over $L_2$
admits a solution is decidable and belongs to $\NP$.
Moreover, the following holds:
\begin{itemize}
\item
If $\sigma_x(w)\ne 0$, then the problem is decidable 
in $\wO(|w|^2)$ time.
\item
If $\sigma_x(w)=0$ and $r(z,y)=0$, then the problem is decidable 
in $2^{O(|w|)}$ time.
\item
If $\sigma_x(w)=0$ and $r(z,y)\ne 0$, then the problem is decidable
in $\wO(|w|^5)$ time.
\end{itemize}
\end{theorem}

\begin{proof}
If $\sigma_x(w)\ne 0$, then decidability in $\wO(|w|^2)$ time
follows from Proposition~\ref{Prop:algorithm-equations}.

Assume now that $\sigma_x(w)=0$.
Compute $\den(w)$ and $\num(w)$.
By Proposition~\ref{pr:compute-num-den}, this can be done
in $O(|w|^2)$ time.
Normalize $\den(w)$ and $\num(w)$, and denote the resulting
$\delta$-parametric polynomials by $\blf$ and $\blg$ respectively.
Compute and normalize $\rev(\blf)$ and $\rev(\blg)$.
Denote the resulting $\delta$-parametric polynomials by
$\blf'$ and $\blg'$ respectively.
Construct the corresponding normalized bivariate polynomials
$F,G,F',G'$ associated with $\blf$, $\blg$, $\blf'$, and $\blg'$, respectively.

By Proposition \ref{prop:main-equivalence},
the equation $w=1$ admits a solution if and only if
$(\den(w), \num(w))$ is a positive instance of $\DIV$.
By Lemma~\ref{le:Q-R-rational}(c),
the remainder $r$ of division of $G$ by $F$ can be computed
within time given by \eqref{eq:u_ij-complexity}.
By Lemma~\ref{le:delta-bound2},
that can be estimated as $\wO(|w|^4)$
in terms of the length of $w$.
By Lemma~\ref{le:r=r'=0},
the remainder $r'$ obtained by dividing $G'$ by $F'$
is zero if and only if $r$ is zero.
Hence, $r'$ need not be computed.

Next, using Proposition~\ref{pr:div-div2}, we compute
$\blf_0\mid\blg_0$,
$\DIV_+(\blf,\blg)$ and $\DIV_+(\rev(\blf),\rev(\blg))$.
By Lemma~\ref{le:rev-characteristics}, 
the instances $(\blf,\blg)$ and $(\rev(\blf),\rev(\blg))$
of $\DIV_+$ have the same characteristics.
By Theorem~\ref{th:DIV-NP},
the problems
$\DIV_+(\blf,\blg)$ and
$\DIV_+(\rev(\blf),\rev(\blg))$
belong to $\NP$.

\medskip\noindent\textsc{Case 1.}
If $r(z,y)=0$, then by Proposition~\ref{pr:r=0-complexity},
$\DIV_+(\blf,\blg)$ and $\DIV_+(\rev(\blf),\rev(\blg))$ can be computed 
in time $2^{O(\deg(\cont(F)))}$, which can be estimated as 
$2^{O(|w|)}$ by Lemma~\ref{le:delta-bound2}.

\medskip\noindent\textsc{Case 2.}
If $r(z,y)\ne 0$, then by Proposition~\ref{pr:r-not0-complexity},
$\DIV_+(\blf,\blg)$ and $\DIV_+(\rev(\blf),\rev(\blg))$ can be computed 
in time \eqref{eq:r-not0-complexity}, which 
can be estimated as $\wO(|w|^5)$ by Lemma~\ref{le:delta-bound2}.
\end{proof}

\section{Generic-case complexity of the Diophantine problem for one-variable equations over $\MZ_2\wr\MZ$}
\label{se:generic-complexity}

In this section we argue that a generic (typical) equation $w(a,t,x) = 1$
satisfies $\sigma_x(w)\ne 0$ and, hence, can be solved in nearly quadratic time.

Let us review several definitions related to generic-case complexity
\cite[Chapter 10]{MSU_book:2011}.
For $m\in\MN$ define the (finite) set
$$
S_{m} = \Set{w\in F(a,t,x)}{|w|=m},
$$
called the \emph{sphere} of radius $m$ in the Cayley graph of $F(a,t,x)$.
Let $P_m$ be the uniform probability distribution on $S_m$.
We say that a set $S\subseteq F(a,t,x)$ is generic in $F(a,t,x)$ if
$$
P_m(S\cap S_m) = \frac{|S\cap S_m|}{|S_m| } \to 1 \mbox{ as } m\to\infty,
$$
i.e., if the density of $S$ within the sets $S_m$
converges to $1$.
We say that a property $P$ of words $w\in F(a,t,x)$ is \emph{generic} if the set 
$$
S=\Set{w}{P \mbox{ holds for } w}
$$
is generic. We say that a \emph{typical word} satisfies $P$
if $P$ is generic.

\subsection{Random one-variable equations $w(a,t,x)$}

A uniformly random word $w(a,t,x) \in S_m$ of length $m$
can be generated as a sequence $w_1,\dots,w_m$ of letters
from the alphabet $\{a^{\pm}, t^{\pm}, x^{\pm}\}$ with no cancellations,
according to the following procedure:
\begin{itemize}
\item
Choose $w_1$ uniformly at random from $\{a^{\pm}, t^{\pm}, x^{\pm}\}$.
\item
For each $i \ge 1$, choose $w_{i+1}$ uniformly at random from
$\{a^{\pm}, t^{\pm}, x^{\pm}\} \setminus \{w_i^{-1}\}$.
\end{itemize}
This generation procedure can be modeled by a discrete-time Markov chain
$\Gamma$ with state set
\[
Q=\{\varepsilon, a^{\pm}, t^{\pm}, x^{\pm}\},
\]
initial state $\varepsilon$, and transition probabilities defined as follows:
\begin{itemize}
\item
From the initial state $\varepsilon$, the chain moves to any state
$q \in \{a^{\pm}, t^{\pm}, x^{\pm}\}$ with uniform probability $\tfrac{1}{6}$.
\item
From any state $q_i \in \{a^{\pm}, t^{\pm}, x^{\pm}\}$, the chain moves to a state
$q_{i+1} \in \{a^{\pm}, t^{\pm}, x^{\pm}\} \setminus \{q_i^{-1}\}$
with uniform probability $\tfrac{1}{5}$.
\end{itemize}
The Markov chain $\Gamma$ thus generates a sequence of random states
$q_0,q_1,q_2,\dots \in Q$, where $q_0=\varepsilon$ almost surely.
For each $m\in\MN$, the word
\[
w = q_1 q_2 \cdots q_m
\]
is a uniformly distributed reduced word of length $m$ in $S_m$.
Our goal is to analyze the random variable $\sigma_x(w)$ for words $w$
generated by this process.

Observe that
\begin{equation}\label{eq:sigma_x}
\sigma_x(q_1\cdots q_m)
=
\bigl|\{1 \le i \le m \mid q_i = x\}\bigr|
-
\bigl|\{1 \le i \le m \mid q_i = x^{-1}\}\bigr|.
\end{equation}
Since the Markov chain $\Gamma$ possesses a high degree of symmetry and
the states $a^{\pm}$ and $t^{\pm}$ do not contribute to the right-hand side
of~\eqref{eq:sigma_x}, we may identify these four states with a single state,
denoted by $y$.
This yields a reduced transition diagram, denoted by $\Gamma^{\ast}$,
shown below.
Equivalently, one may view the resulting process as generating words
in which each occurrence of $a^{\pm}$ or $t^{\pm}$ is replaced by the
single letter $y$, while occurrences of $x^{\pm}$ are left unchanged.

\begin{figure}[h]
\centering
\begin{minipage}[b]{0.48\linewidth}
\centering
\includegraphics[width=\linewidth]{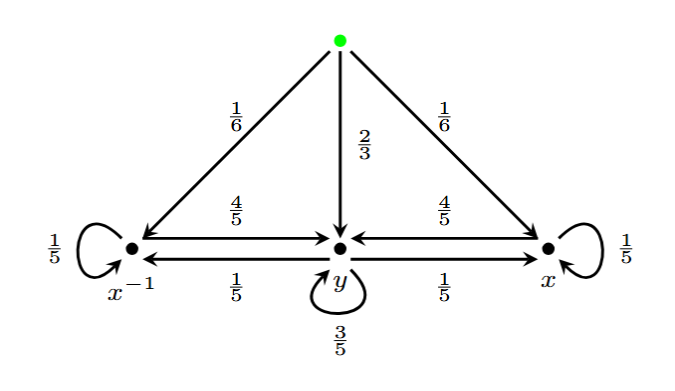}
\captionof{subfigure}{Transition diagram $\Gamma^\ast$.}
\label{fig:transition_probability_with_epsilon}
\end{minipage}
\hfill
\begin{minipage}[b]{0.48\linewidth}
\centering
\includegraphics[width=\linewidth]{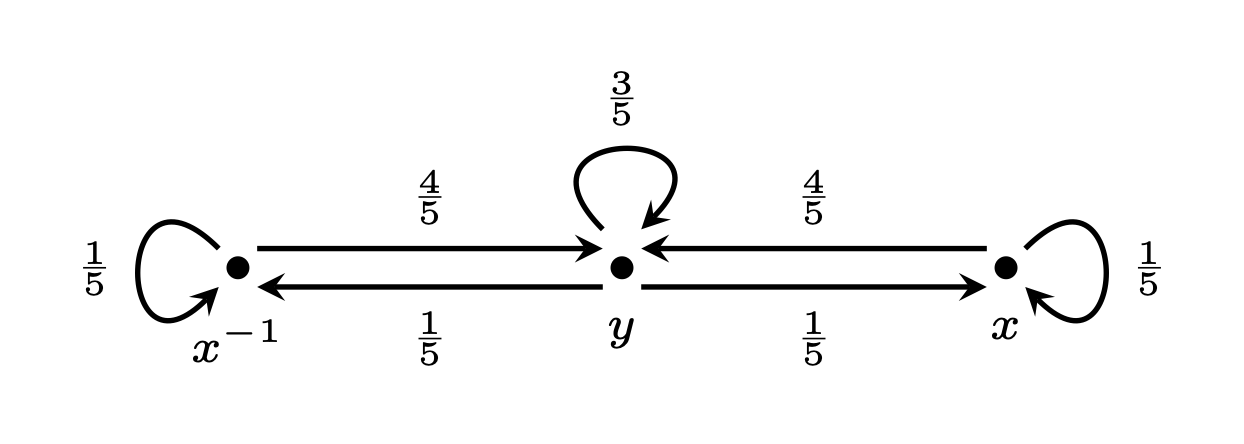}\vspace{5mm}
\captionof{subfigure}{Irreducible $\Gamma^\ast$.}
\label{fig:transition_probability}
\end{minipage}
\end{figure}

Notice that this process is not irreducible. However, after deleting
the transient state $\varepsilon$, we obtain a new Markov chain
$q_1', q_2', \dots$ whose initial distribution is given by
\[
\Pr[q_1' = q] =
\begin{cases}
\frac{1}{6}, & \text{if } q = x^{-1},\\[2mm]
\frac{2}{3}, & \text{if } q = y,\\[2mm]
\frac{1}{6}, & \text{if } q = x.
\end{cases}
\]
The transition matrix of this chain is
\[
M =
\begin{bmatrix}
0.2 & 0.8 & 0 \\
0.2 & 0.6 & 0.2 \\
0   & 0.8 & 0.2
\end{bmatrix}.
\]
The resulting Markov chain is \emph{irreducible}, since every state
can be reached from any other state with positive probability, and
\emph{aperiodic}, because the greatest common divisor of all return
times to any state equals~$1$. Consequently, the chain is
\emph{ergodic}. In particular,
\begin{itemize}
\item
it admits a unique \emph{stationary distribution} $\pi$
(a probability vector satisfying $\pi M = \pi$), which in this case is
\[
\pi = \Bigl(\tfrac{1}{6}, \tfrac{2}{3}, \tfrac{1}{6}\Bigr),
\]
\item
and, starting from any initial distribution, the distribution of
$q_n'$ converges to $\pi$ as $n \to \infty$.
\end{itemize}
It is well known that ergodic Markov chains on finite state spaces are
\emph{uniformly ergodic}, meaning that convergence to the stationary
distribution occurs at a uniform geometric rate, independent of the
initial state. Moreover, the chain is \emph{Harris recurrent}: starting
from any initial state, it visits every set of positive stationary
measure infinitely often with probability~$1$.

Now consider the function $f \colon \{x^{-1}, y, x\} \to \{-1,0,1\}$ defined by
\[
f(x^{-1}) = -1, \qquad f(y) = 0, \qquad f(x) = 1.
\]
Clearly, $\ME_{\pi}[f] = 0$. By construction of $f$, we have
\begin{equation}\label{eq:sigma_x2}
\sigma_x(q_1' \dots q_m') = \sum_{i=1}^{m} f(q_i').
\end{equation}
Define the random variable
\[
\bar{f}_m = \frac{1}{m} \sum_{i=1}^{m} f(q_i'),
\]
which is the empirical (time-averaged) estimator of $\ME_{\pi}[f]$ along the trajectory
$q_1',\dots,q_m'$.

\begin{theorem*}[Central limit theorem for Markov chains {\cite[Theorem~9(6)]{Jones:2004}}]
Let $X = (X_n)_{n\ge 1}$ be a Harris ergodic Markov chain with stationary
distribution $\pi$, and let $f$ be a real-valued function satisfying
$\ME_\pi[f^2] < \infty$.  
If $X$ is uniformly ergodic, then
\[
\sqrt{m}\,\bigl(\bar f_m - \ME_\pi[f]\bigr)
\;\xrightarrow{\,d\,}\;
\mathcal N(0,\sigma^2),
\]
for some $\sigma^2 \ge 0$.
\end{theorem*}

For our Markov chain, the central limit theorem gives
\[
\frac{1}{\sqrt{m}} \sum_{i=1}^m f(q_i') 
\;\stackrel{d}{\longrightarrow}\; \mathcal{N}(0,\sigma^2),
\]
which immediately implies the following.

\begin{proposition}
$\Pr[f(q_1')+\dots+f(q_m') = 0]\to 0$ as $m\to \infty$.
\end{proposition}

\begin{proof}
We have
\[
\Pr\Big[\sum_{i=1}^m f(q_i') = 0\Big] 
= \Pr\Bigg[\frac{1}{\sqrt{m}} \sum_{i=1}^m f(q_i') = 0\Bigg].
\]
By the central limit theorem, the distribution of 
$\frac{1}{\sqrt{m}} \sum_{i=1}^m f(q_i')$ converges to the normal distribution 
$\mathcal{N}(0,\sigma^2)$. Since the probability that a continuous normal 
variable equals exactly zero is $0$, the result follows.
\end{proof}

\begin{corollary}\label{co:sigma-x-generic}
$P_m[\Set{w\in S_m}{\sigma_x(w)=0}] \to 0 \mbox{ as } m\to \infty.$
\end{corollary}

\begin{corollary}
The problem
$\DP_1(L_2)$ has $O(|w|^2 \log |w| \log \log |w|)$ generic-case time complexity.
\end{corollary}

\begin{proof}
By Proposition~\ref{Prop:algorithm-equations}, there exists an algorithm that, for any equation $w = 1$ with $\sigma_x(w) \neq 0$, determines in 
$O(|w|^2 \log |w| \log \log |w|)$ time whether a solution exists.
Moreover, by Corollary~\ref{co:sigma-x-generic}, the set of words
\[
\{ w \in F(a,t,x) \mid \sigma_x(w) \neq 0 \}
\]
is generic in $F(a,t,x)$.
\end{proof}

\bibliography{main_bibliography}

\end{document}